\documentclass[aos]{imsart}

\RequirePackage{amsthm,amsmath,amsfonts,amssymb,tabularx}
\RequirePackage[numbers,sort&compress]{natbib}
\RequirePackage[colorlinks,citecolor=blue,urlcolor=blue]{hyperref}
\RequirePackage{graphicx,bbm,bm}
\RequirePackage[graphicx]{realboxes}
\RequirePackage{caption,subcaption}
\RequirePackage{float}
\RequirePackage{xr-hyper}
\startlocaldefs

\makeatletter

\theoremstyle{plain}

\newtheorem{theorem}{Theorem}[section]

\newtheorem{proposition}{Proposition}
\theoremstyle{remark}
\newtheorem{assumption}{Assumption}[section]
\newtheorem{definition}[theorem]{Definition}
\newtheorem{example}{Example}

\newtheorem{remark}{Remark}[section]
\pagebreak[0]
\endlocaldefs

\begin{document}

\begin{frontmatter}
\title{Generalized Linear Spectral Statistics of High-dimensional Sample Covariance Matrices and Its Applications}
\runtitle{GLSS of High-dimensional Sample Covariance Matrices}

\begin{aug}
\author[A]{\fnms{Yanlin}~\snm{Hu}\ead[label=e1]{hyl11@mail.ustc.edu.cn}},
\author[A]{\fnms{Qing}~\snm{Yang}\ead[label=e2]{ yangq@ustc.edu.cn}}
\and
\author[A]{\fnms{Xiao}~\snm{Han}\ead[label=e3]{xhan011@ustc.edu.cn}}
\address[A]{International Institute of Finance, School of Management, University of Science and Technology of China\printead[presep={,\ }]{e1,e2,e3}}

\end{aug}

\begin{abstract}
		In this paper, we introduce the \textbf{G}eneralized \textbf{L}inear \textbf{S}pectral \textbf{S}tatistics (GLSS) of a high-dimensional sample covariance matrix $\bm{S}_n$, denoted as $\operatorname{tr}f(\bm{S}_n)\bm{B}_n$, which effectively captures distinct spectral properties of $\bm{S}_n$ by incorporating an ancillary matrix $\bm{B}_n$ and a test function $f$.
		The joint asymptotic normality of GLSS associated with different test functions is established under mild assumptions on $\bm{B}_n$ and the underlying distribution, when the dimension $n$ and sample size $N$ are comparable.  The convergence rate of GLSS is determined by $\sqrt{{N}/{\operatorname{rank}(\bm{B}_n)}}$. Subsequently, we propose a novel functional projection approach based on GLSS for hypothesis testing on eigenspaces of ``population-spiked'' covariance matrices, showcasing a universality phenomenon in the magnitude of the spikes. The theoretical accuracy of our results established for GLSS and the advantages of the newly suggested testing procedure are
demonstrated through various numerical studies. 		
\end{abstract}

\begin{keyword}[class=MSC]
\kwd[Primary ]{62H10}
\kwd{60B20}
\kwd[; secondary ]{62H15}
\kwd{60F05}
\end{keyword}

\begin{keyword}
\kwd{Sample covariance matrix}
\kwd{random matrix theory}
\kwd{eigenspaces}
\kwd{generalized linear spectral statistics}
\end{keyword}

\end{frontmatter}

\section{Introduction}\label{sec1}
The covariance matrix holds paramount importance in statistics and its associated fields, serving as a fundamental component for numerous widely-used methodologies that heavily rely on comprehending its structural characteristics. For instance, methodologies such as principal component analysis \citep{Pearson:1901} and factor analysis \citep{Fan:2011, Li:2018} depend on understanding the eigenstructure corresponding to the leading eigenvalues, while spectral methods in clustering \citep{han2023eigen} depend on understanding the asymptotic properties of the eigenvectors containing the clustering information. Although the sample covariance matrix is a consistent estimator of its population counterpart in the low-dimensional setting with a fixed number of variables $n$, it is widely recognized that drawing direct inferences from the sample covariance matrix may lead to erroneous conclusions when the dimensionality $n$ is comparable to or significantly larger than the sample size $N$ \citep{Yao:2015}. Specifically, for example, \citep{Jone:2001,Jone:2009,Wachter:1978} have shown that when
	$n/N\rightarrow c\in (0,\infty)$, the largest eigenvalue of the sample covariance matrix is an inconsistent estimator for the largest eigenvalue of the population covariance matrix, and the eigenvectors
	of the sample covariance matrix can be nearly orthogonal to the true ones.

In the high-dimensional setting, numerous monographs have been dedicated to investigating the asymptotic behavior of the largest few eigenvalues or the spectrum of sample covariance matrices. \citep{BY:1993} and \citep{Yin:1988} established the almost sure convergence to the edge of Marchenko-Pastur (M-P) law for the {smallest and largest} eigenvalues of sample covariance matrices, respectively. Subsequently, many efforts have been devoted to characterizing the asymptotic distribution of the largest eigenvalue or joint distribution of a few leading eigenvalues. We refer the readers to the literatures \citep{Baik:2006,Blo:2016,Cai:2020,Jone:2001,El:2007,Lee:2016,Paul:2007} and the references therein for more detailed discussions. Regarding the spectrum, \citep{Bai:2004} established the central limit theorem (CLT) for linear spectral statistics of sample covariance matrices, which considers the sum of eigenvalues of $f(\bm{S}_n)$ (i.e. $\operatorname{tr}f(\bm{S}_n)$), where $f$ is assumed to be analytic. The Gaussian-like fourth moment assumption therein and the constrains made on the test function $f$ were later relaxed by \citep{Najim:2015,Pan:2008,Zheng:2015}. Many statistical inference problems on population covariance matrices can be addressed by employing the CLT of linear spectral statistics, as exemplified in studies \citep{Bai:2009,Zheng:2019}.

In recent years, there has been a growing interest regarding the properties of eigenvectors of sample covariance matrices. Under different assumptions on the structure of population covariance matrices and on the distribution of underlying variables, various works have focused on deriving the asymptotic behavior of the inner product between eigenvectors of sample covariance matrices and some non-random vectors. To name a few, we refer the readers to  \citep{Blo:2016,Cai:2020,Jone:2009,John:2018,Liu:2023a,Mestre:2008,Paul:2007,Yin:2023}. Recently, \citep{Bao:2022} established the asymptotic expansion of the spiked
	eigenvalues and linear combination of spiked eigenvectors for a high-dimensional spiked
	covariance matrix model. Their theoretical results necessitate that the non-spiked part of the population covariance matrix is an identity matrix, while also assuming a finite number of spiked eigenvalues and arbitrary finite moments for the data entries.
 \citep{Bai:2007} proposed another statistic to analyze eigenvalues and eigenvectors by introducing a non-random unit test vector $\bm{b}_n$.  To be more specific, they conducted an investigation on $\bm{b}_n^*f(\bm{S}_n)\bm{b}_n$ and established its CLT, while referring to \citep{Pan:2008} for a related work under weaker assumptions.

The purpose of the present paper is to establish the CLT for \textbf{G}eneralized \textbf{L}inear \textbf{S}pectral \textbf{S}tatistics (GLSS) of sample covariance matrices, which is formally defined as follows:
	\begin{equation}\label{core}
		\operatorname{tr}f(\bm{S}_n)\bm{B}_n,
	\end{equation}
	where the sample covariance matrix $\bm{S}_n$ takes the form
	\begin{equation}\label{xh1}
		\bm{S}_n=\frac{1}{N}\bm{\Sigma}_n^{1/2}\bm{X}_n\bm{X}_n^*\bm{\Sigma}_n^{1/2}=\frac{1}{N}\sum_{j=1}^N\bm{\Sigma}_n^{1/2}\bm{x}_j\bm{x}_j^*\bm{\Sigma}_n^{1/2}.
	\end{equation}
The entries of the $n\times N$ matrix $\bm{X}_n=(X_{i,j}^n)$ are i.i.d with zero mean and unit variance and $\bm{x}_j=(X_{1,j}^n,\ldots,X_{n,j}^n)^\top$, $j=1,\ldots,N$. The matrix $\bm{\Sigma}_n^{1/2}$ represents the square root of the population covariance matrix $\bm{\Sigma}_n$.
When $\bm{B}_n$ equals to the identity matrix $\bm{I}_n$, GLSS is the standard linear spectral statistics introduced in  \citep{Bai:2004}. In the case of $\bm{B}_n$ being a rank one Hermitian matrix, GLSS reduces to the statistic considered in \citep{Bai:2007}.
We mention {five} other relevant works. Firstly, \citep{Cipo:2023} established the CLT for $\operatorname{tr}f(\bm{W}_n)\bm{B}_n$, where $\bm{W}_n$ is a Wigner matrix. Given the existence of arbitrary finite moments and $\|\bm{B}_n\|_F\geq cn^{\epsilon}$ for some $c,\epsilon>0$, they proved that $\operatorname{tr}f(\bm{W}_n)\bm{B}_n$ is asymptotic Gaussian. While in our paper, we will develop the CLT of GLSS by considering the existence of the fourth moment and exploring various ranks of $\bm{B}_n$ under mild assumptions on its structure. Secondly, \citep{Le:2010} determined the almost sure limit of $\operatorname{tr}\left(\bm{S}_n-z \bm{I}_n\right)^{-1}g\left(\bm{\Sigma}_n\right)$ for some bounded function $g$ and complex number $z$, under the condition $\mathbb{E}|X_{ij}|^{12}<\infty$. {Thirdly, \citep{Rubio:2011} considered a general class of random matrices taking the form: $\bm{S}_{n,ge}=\bm{A}_n+N^{-1}\bm{\Sigma}_n^{1/2}\bm{X}_n\bm{T_n}\bm{X}_n^*\bm{\Sigma}_n^{1/2}$ where $\bm{A}_n$, $\bm{\Sigma}_n$ and
$\bm{T}_n$ are Hermitian nonnegative definite matrices, such that $\bm{\Sigma}_n$ and $\bm{T}_n$ have bounded spectral
norm with $\bm{T}_n$ being diagonal. They determined the almost sure limit of $\operatorname{tr}\left(\bm{S}_{n,ge}-z \bm{I}_n\right)^{-1}\bm{B}_n$ by assuming $\mathbb{E}|X_{ij}|^{8+\varepsilon}<\infty$ for some $\varepsilon>0$ and $\|\bm{B}_n\|_F<\infty$. Fourthly, \citep{Bodnar:2024} derived the almost sure limit of weighted moments of Moore-Penrose inverse and the ridge-type inverse of the centered sample covariance matrix.}  Lastly, in our parallel working paper, we have developed the CLT of GLSS for high-dimensional sample correlation matrices. {To our best knowledge, this is the first
work concerning the asymptotic distribution of $\operatorname{tr}f(\bm{S}_n)\bm{B}_n$ for general $\bm{B}_n$. }

Given that the matrix $\bm{B}_n$ has a rank of $k_n$ and possesses a spectral decomposition $\bm{B}_n=\sum_{i=1}^{k_n}s_i\bm{b}_i\bm{b}_i^*$, GLSS (\ref{core}) can be rewritten as
	\begin{equation}\label{eg1}
\sum_{i=1}^{k_n}s_i\bm{b}_i^*f(\bm{S}_n)\bm{b}_i,
\end{equation}
	which is a weighted sum of the vector linear spectral statistics (i.e. $k_n=1$) considered by \citep{Bai:2007,Pan:2008}.  The extension to general $k_n$, especially when $k_n$ diverges with $n$ is non-trivial and the proof is much more complicated. A further spectral decomposition on $\bm{S}_n$ yields an alternative representation of GLSS as follows:
	\begin{equation}\label{sep}
		\sum_{i=1}^{k_n}\sum_{j=1}^n s_i f(\lambda_j)|\langle\bm{b}_i,\bm{u}_j\rangle|^2,
	\end{equation}
	where $\lambda_j$ is the $j$-th largest eigenvalue of $\bm{S}_n$ and $\bm{u}_j$ is the corresponding eigenvector. It is evident from \eqref{sep} that by selecting different choices of $f$ and $\bm{B}_n$, GLSS reflects distinct aspects of the spectrum of $\bm{S}_n$. Therefore it becomes feasible to assess a partial spectral structure of $\bm{S}_n$ through appropriate selection of $f$ and $\bm{B}_n$. This has been verified in our application, where we propose a novel approach - functional projection - for conducting hypothesis testing on eigenspaces of ``population-spiked'' covariance matrices. The concept of population-spike will be elaborated upon extensively therein.

	Our main contributions can be summarized as follows:
\begin{itemize}
		\item[$\bullet$]We propose a flexible statistic - GLSS to
  study the properties of eigenvalues and eigenvectors of high-dimensional sample covariance matrices. The statistics studied in \citep{Bai:2007} and \citep{Bai:2004} are special cases of GLSS.
		
		\item[$\bullet$] We establish the CLT of GLSS for all $1\leq k_n\leq n$ using an adaptive proof procedure for different values of $k_n$. Notably, we introduce a new two-step truncation strategy when dealing with the case where $k_n/n\rightarrow 0$. Moreover, we relax the assumptions made in  \citep{Bai:2007} and \citep{Pan:2008}, which specifically consider the case where $k_n=1$; please refer to Remark \ref{xrmk2} for further details. Due to the existence of $f$ and $\bm{B}_n$, this CLT helps to understand the
  eigenvalue and eigenvector structure of $\bm{S}_n$ in a flexible way. For instance, choosing $\bm{B}_n$ as a projection matrix allows $f(\bm{S}_n)\bm{B}_n$ to represent the projection of $f(\bm{S}_n)$ onto the space of $\bm{B}_n$. Additionally, by utilizing different ranks $k_n$, we could keep arbitrary number of projection directions.




		\item[$\bullet$] Based on GLSS, we propose a novel functional projection approach for conducting eigenspace testing on covariance matrices with ``population-spiked'' characteristics. The term ``population-spiked'' is employed here to distinguish our method from existing approaches that impose lower bound constraints on the spikes; in contrast, our method accommodates varying numbers of spikes without making such assumptions on their magnitudes.


	\end{itemize}

	The {remainder} of this article is organized as follows. In Section \ref{section2}, we establish the CLT for GLSS by considering both cases when $k_n$ is comparable to $n$ and when $k_n=o(n)$. Various simulations are conducted in Section \ref{section4} to verify our theoretical results.  Motivated by GLSS and building upon a slight modification to our Theorem \ref{main}, we propose a novel test statistic in Section \ref{section5} for testing eigenspaces of population-spiked covariance matrices. To demonstrate the advantages of our proposed method, we conduct comprehensive comparisons with various alternative methods across multiple aspects, including computational complexity, accuracy under null hypotheses, and power under alternative hypotheses.
	All auxiliary lemmas and proofs, as well as additional results are postponed to the supplementary material.

\textbf{Notations.} We introduce some notations that will be used throughout this paper. Bold capital and lowercase letters are used to denote matrices and vectors, respectively.  The notation $\stackrel{D}{\rightarrow}$ (or $\stackrel{\mathbb{P}}{\rightarrow}$) means convergence in distribution (or in probability). {For any quantities $a_n$ and $b_n$, we use the notation $a_n\ll b_n$ to denote the relation $a_n/b_n\rightarrow0$ as $n\rightarrow\infty$. In addition, we write $a_n\asymp b_n$ if there exists some constant $C>1$ such that $C^{-1}|a_n|\leq |b_n|\leq C|a_n|$. For random variable sequences $x_n$, the symbol $x_n=\mathrm{o}_{\mathbb{P}}(a_n)$ means $x_n/a_n\stackrel{\mathbb{P}}{\rightarrow}0$. Besides, $x_n=\mathrm{O}_{\mathbb{P}}(a
_n)$ stands for $\lim_{M\rightarrow\infty}\sup_{n}\mathbb{P}(|x_n/a_n|>M)=0$.} For a matrix $\bm{M}\in \mathbb{C}^{p\times q}$, we use $\|\bm{M}\|$ and $\|\bm{M}\|_F$ to denote its spectral norm and Frobenius norm.  In addition, denote by $(\bm{M})_{ij}$, $\lambda_{i}(\bm{M})$ and $s_i(\bm{M})$ the entry located in the $i$-th row and $j$-th column, the $i$-th largest eigenvalue and the $i$-th largest singular value, respectively.   Let $\bm{M}^*$ (or $\bm{M}^{\top}$) represent the conventional conjugate transpose (or transpose) of $\bm{M}$. The notation $\operatorname{diag}(\bm{M})$  in the context of a square matrix $\bm{M}$  denotes a diagonal matrix whose entries on the main diagonal correspond to those of $\bm{M}$. {For two matrices $\bm{M}$ and $\bm{N}$ of the same size, we write $\bm{M}\circ\bm{N}$ for their Hadamard product.}
For a  $\sigma$-field $ \mathcal{F}_i$ generated by $\{\textbf{x}_1,...,\textbf{x}_i\}$, we use $\mathbb{E}_{i}(\cdot)$ to denote the conditional expectation with respect to $\mathcal{F}_i$.  Furthermore, denote by $\mathbb{I}_{E}$ the indicator function of an event $E$. For a compact metric space $(\mathcal{K},d)$, let $C(\mathcal{K},\|\cdot\|_{\infty})$ represent the space of
			continuous complex-valued functions on $\mathcal{K}$ equipped with the uniform norm, i.e. $\|f\|_{\infty}=\sup_{t\in\mathcal{K}}|f(t)|$.

\section{Asymptotic Results for GLSS}\label{section2}
	In this section, the asymptotic distribution of our GLSS is established both when $\frac{k_n}{n}\rightarrow 0$ and $\frac{k_n}{n}\geq c_0$ for some positive constant $c_0$.  Before delving into the main theorems in Section \ref{sec2.2}, we provide an introduction to some preliminary results regarding the limiting spectral distribution of the conventional sample covariance matrix $\bm{S}_n$ in Section \ref{sec2.1}.

\subsection{{Some preliminary results on the sample covariance matrix}}\label{sec2.1}
	In random matrix theory, the Stieltjes transform is a fundamental function, which is formally defined in Definition \ref{def}.
	\begin{definition}\label{def}
		For any function $G$ with bounded variation on the real line, its Stieltjes transform is defined by
		$$m_G(z)=\int \frac{1}{x-z}dG(x), \qquad z\in\mathbb{C} ~\text{and}~ \Im{z}\neq0.$$
	\end{definition}
It has been demonstrated that a bijective correspondence exists between $G$ and its Stieltjes transform $m_G(z)$ when $G$ is a proper distribution function (see Theorem B.8 in \citep{Bai:2010}).
	Recalling the definition of $\bm{S}_n$ in \eqref{xh1}, an elementary limit theorem concerning the eigenvalues of $\bm{S}_n$ focuses on its empirical spectral distribution $F^{\bm{S}_n}$, which is defined as
	$$
	F^{\bm{S}_n}(x)=\frac{1}{n}\sum_{i=1}^n\mathbb{I}_{\{\lambda_i(\bm{S}_n)\leq x\}}.
	$$
To be more specific, if we assume that for all $n$, $X_{ij}^n$ are i.i.d. random variables with zero mean and unit variance, $H_n=F^{\bm{\Sigma}_n}$ convergences in distribution to $H$, a proper cumulative distribution function (c.d.f.) and $c_n=n/N\rightarrow c\in (0,\infty)$, then almost surely, $F^{\bm{S}_n}$ converges in distribution to $F^{c, H}$, a nonrandom proper c.d.f whose  Stieltjes transform $m(z)$ is the unique solution to
	\begin{equation}\label{m3}
		m(z)=\int \frac{1}{x(1-c-c z m(z))-z} d H(x), \quad z \in \mathbb{C}^{+}.
	\end{equation}
Considering $\bm{\underline{S}}_n \equiv(1 / N) \bm{X}_n^* \bm{\Sigma}_n \bm{X}_n$ whose spectra differs from that of $\bm{S}_n$ by $|n-N|$ zeros, we know that its limiting empirical distribution function satisfies
	$$
	\underline{F}^{c, H} \equiv(1-c) \mathbb{I}_{[0, \infty)}+c F^{c, H}.
	$$
Furthermore, its Stieltjes transform
	\begin{equation}\label{m1}
		\underline{m}(z) \equiv m_{\underline{F}^{c, H}}(z)=-\frac{1-c}{z}+c m(z)
	\end{equation}
	has inverse
	\begin{equation}\label{m2}
		z=z(\underline{m})=-\frac{1}{\underline{m}}+c \int \frac{t}{1+t \underline{m}} d H(t),
	\end{equation}
	which takes a simpler form. One may refer to \citep{Bai:2010} for more detailed discussions.
	Let ${m}_n^0(z)$ and $\underline{m}_n^0(z)$ represent the quantities obtained from equations (\ref{m3}) and (\ref{m1}) when replacing $(c, H)$ by $(c_n, H_n)$, which will be frequently used in establishing our main theorems. The corresponding distribution functions for $m_n^0(z)$ and $\underline{m}_n^0(z)$ are denoted as $F^{c_n,H_n}$ and $\underline{F}^{c_n,H_n}$, respectively. In addition, $m_n(z)$ and $\underline{m}_n(z)$ are employed to denote the Stieltjes transforms of $F^{\bm{S}_n}$ and $F^{\underline{\bm{S}}_n}$.

	\subsection{{Main theoretical results}}\label{sec2.2}
The following assumptions will be used in our theoretical analysis.
\begin{assumption}\label{asa}
{For each $n$, $ X_{i j}=X_{ij}^n,\ 1\leq i \leq n,\ 1\leq j \leq N$, are i.i.d. for all $i, j$.} Moreover, $\mathbb{E} X_{11}=0$, $\mathbb{E}\left|X_{11}\right|^2=1$, $\mathbb{E}\left|X_{11}\right|^4<\infty$, $c_n=n / N \rightarrow c\in(0,\infty)$. For complex case we assume $\mathbb{E}X_{11}^2=0$.
\end{assumption}

 \begin{assumption}\label{asb}
 The matrices $\bm{\Sigma}_n$ and $\bm{B}_n$ are $n\times n$ non-random Hermitian matrices such that their non-zero eigenvalues are bounded away from $0$ and infinity. Moreover, we assume that $\bm{\Sigma}_n$ is non-negative definite ($\bm{\Sigma}_n\succeq \bm{0}$) and  $H_n=F^{\bm{\Sigma}_n} \stackrel{D}{\rightarrow} H$, where $H$ is a proper c.d.f.
 \end{assumption}

  \begin{assumption}\label{asc}
  Let $k_n=\operatorname{rank}(\bm{B}_n)$. Either one of the following two cases holds:

   (i) {[$k_n$ is comparable to $n$].} There exists a positive constant $c_0$ such that
		$\frac{k_n}{n}\geq c_0.$

  (ii) {[$k_n$ is much smaller than $n$].} $\frac{k_n}{n}\rightarrow  0.$
 \end{assumption}

Assumptions \ref{asa} and \ref{asb} are standard in random matrix theory (see \citep{Bai:2007,Bai:2004,Pan:2008} for example). The asymptotic behavior of GLSS -- $\operatorname{tr}f(\bm{S}_n)\bm{B}_n$ depends on the rank of $\bm{B}_n$. In the following Theorem \ref{main} and Theorem \ref{the6.1}, we summarize the different limiting distributions under the two different cases of $k_n$ stated in Assumption \ref{asc}, respectively.
{Before presenting the main results, we introduce the following definitions used therein. Define $\bm{\overline{\Sigma}}_n(z)=\bm{I}_n+\underline{m}_n^0(z)\bm{\Sigma}_n$ and the following quantities:
\begin{align}
    P_n(z)&=\frac1N\operatorname{tr}\left(\bm{\overline{\Sigma}}_n^{-2}(z)\bm{\Sigma}_n\bm{B}_n\right),\qquad
                Q_n(z)=\frac1N\operatorname{tr}\left(\bm{\overline{\Sigma}}_n^{-3}(z)\bm{\Sigma}_n^2\bm{B}_n\right),\notag\\V_n^3(z_1,z_2)&=\frac1{z_1^2z_2^2N}\operatorname{tr}\left(\bm{\overline{\Sigma}}_n^{-1}(z_2)\bm{\overline{\Sigma}}_n^{-1}(z_1)\bm{\Sigma}_n\bm{B}_n\bm{\overline{\Sigma}}_n^{-1}(z_1)\bm{\overline{\Sigma}}_n^{-1}(z_2)\bm{\Sigma}_n\bm{B}_n\right),\notag\\\widetilde{V}_n^1(z_1,z_2)&=\frac{1}{z_1z_2^2N}\sum_{i=1}^n\left(\bm{\overline{\Sigma}}_n^{-1}(z_1)\bm{\Sigma}_n\right)_{ii}\left(\bm{\Sigma}_n^{1/2}\bm{\overline{\Sigma}}_n^{-1}(z_2)\bm{B}_n\bm{\overline{\Sigma}}_n^{-1}(z_2)\bm{\Sigma}_n^{1/2}\right)_{ii},\notag\\\widetilde{V}_n^2(z_1,z_2)&=\frac{1}{z_1^2z_2^2N}\sum_{i=1}^n\left(\bm{\overline{\Sigma}}_n^{-2}(z_1)\bm{\Sigma}_n^2\right)_{ii}\left(\bm{\Sigma}_n^{1/2}\bm{\overline{\Sigma}}_n^{-1}(z_2)\bm{B}_n\bm{\overline{\Sigma}}_n^{-1}(z_2)\bm{\Sigma}_n^{1/2}\right)_{ii},\notag \\V_n^1(z_1,z_2)&=\frac1{z_1z_2^2N}\operatorname{tr}\left(\bm{\overline{\Sigma}}_n^{-2}(z_2)\bm{\overline{\Sigma}}_n^{-1}(z_1)\bm{\Sigma}_n^2\bm{B}_n\right),\notag\\V_n^2(z_1,z_2)&=\frac1{z_1^2z_2^2N}\operatorname{tr}\left(\bm{\overline{\Sigma}}_n^{-2}(z_2)\bm{\overline{\Sigma}}_n^{-2}(z_1)\bm{\Sigma}_n^3\bm{B}_n\right),\notag\\                U_n^1(z_1,z_2)&=\frac1{z_1z_2^2N}\operatorname{tr}\left(\bm{\overline{\Sigma}}_n^{-2}(z_2)\bm{\overline{\Sigma}}_n^{-1}(z_1)\bm{\Sigma}_n^3\right),\notag\\
U_n^2(z_1,z_2)&=\frac1{z_1^2z_2^2N}\operatorname{tr}\left(\bm{\overline{\Sigma}}_n^{-2}(z_2)\bm{\overline{\Sigma}}_n^{-2}(z_1)\bm{\Sigma}_n^4\right),\label{za}\\g_n(z)&=\frac{P_{n}(z)}{z^2}\left(1-\frac{(\underline{m}_n^0(z))^2}{N}\operatorname{tr}\bm{\overline{\Sigma}}_n^{-2}(z)\bm{{\Sigma}}_n^2\right)^{-1},\notag\\a_n(z_1,z_2)&=\frac{\underline{m}_n^0(z_1)\underline{m}_n^0(z_2)}{N}\operatorname{tr}\left(\bm{\overline{\Sigma}}_n^{-1}(z_1)\bm{\overline{\Sigma}}_n^{-1}(z_2)\bm{\Sigma}_n^2\right),\notag\\\widetilde{U}_n^1(z_1,z_2)&=\frac{1}{z_1z_2^2N}\sum_{i=1}^n\left(\bm{\overline{\Sigma}}_n^{-1}(z_1)\bm{\Sigma}_n\right)_{ii}\left(\bm{\overline{\Sigma}}_n^{-2}(z_2)\bm{\Sigma}_n^2\right)_{ii},\notag\\\widetilde{U}_n^2(z_1,z_2)&=\frac{1}{z_1^2z_2^2N}\sum_{i=1}^n\left(\bm{\overline{\Sigma}}_n^{-2}(z_1)\bm{\Sigma}_n^2\right)_{ii}\left(\bm{\overline{\Sigma}}_n^{-2}(z_2)\bm{\Sigma}_n^2\right)_{ii},\notag\\\widetilde{V}_n^3(z_1,z_2)&=\frac{1}{z_1^2z_2^2N}\sum_{i=1}^n\left(\bm{\Sigma}_n^{1/2}\bm{\overline{\Sigma}}_n^{-1}(z_1)\bm{B}_n\bm{\overline{\Sigma}}_n^{-1}(z_1)\bm{\Sigma}_n^{1/2}\right)_{ii}\notag\\&\qquad\qquad\times\left(\bm{\Sigma}_n^{1/2}\bm{\overline{\Sigma}}_n^{-1}(z_2)\bm{B}_n\bm{\overline{\Sigma}}_n^{-1}(z_2)\bm{\Sigma}_n^{1/2}\right)_{ii},\notag\\\tilde{a}_n(z_1,z_2)&=\frac{\underline{m}_n^0(z_1)\underline{m}_n^0(z_2)}{N}\sum_{i=1}^n\left(\bm{\overline{\Sigma}}_n^{-1}(z_1)\bm{\Sigma}_n\right)_{ii}\left(\bm{\overline{\Sigma}}_n^{-1}(z_2)\bm{\Sigma}_n\right)_{ii},\notag\\\zeta_n^1(z_1,z_2)&=V_n^1(z_1,z_2)+z_2^2\underline{m}_n^0(z_2)g_n(z_2)U_n^1(z_1,z_2).\notag\end{align}

}
	\begin{theorem}\label{main} [$k_n$ is comparable to $n$]. Suppose that  Assumptions \ref{asa}, \ref{asb} and \ref{asc} (i) hold.
	   Let $f_1, \ldots, f_r$ be analytic functions  on an open interval containing $[d_{-},d^{+}]$, where
		\begin{equation}\label{interval}
			[d_{-},d^{+}]=\left[\liminf_{n}\lambda_{\min}^{\bm{\Sigma}_n}\mathbb{I}_{(0,1)}(c)(1-\sqrt{c})^2, \limsup_{n}{\lambda_{\max}^{\bm{\Sigma}_n}}(1+\sqrt{c})^2\right].
		\end{equation}
Recall the definition of GLSS in \eqref{core} and define
\begin{equation}\label{add2-1}
  \Theta_n(f)=\operatorname{tr} f(\bm{S}_n)\bm{B}_n-\frac{1}{2\pi i}\oint_{\Gamma}f(z)\operatorname{tr}(z\bm{I}_n+z\underline{m}_n^0(z)\bm{\Sigma}_n)^{-1}\bm{B}_ndz,
\end{equation}
 where $\Gamma$ is a contour taken in the positive direction enclosing an open interval covering $[d_{-},d^{+}]$.
  Then we have the following results:
		
		\vskip 0.4cm
		(i) the random vector
		\begin{equation}\label{Vector}
			\left(\Theta_n(f_1), \ldots, \Theta_n(f_r)\right)
		\end{equation}
		forms a tight sequence in $n$.
		
		(ii) Let $\mu_{X}=\mathbb{E}|X_{11}|^4-\left|\mathbb{E}X_{11}^2\right|^2-2$ and $\upsilon_{X}=1+\left|\mathbb{E}X_{11}^2\right|^2$.  After suitable centralization, the random vector (\ref{Vector}) converges weakly to an $r$-dimensional Gaussian distribution, i.e.,
		\begin{align}\label{converge}
			(\Theta_n(f_1)-\omega_n(f_1),\cdots,\Theta_n(f_r)-\omega_n(f_r))\stackrel{D}{\rightarrow}\mathcal{N}(\bm{0},\bm{\Omega}_1),
		\end{align}
		where
		\begin{align}\label{meansam}
			\begin{split}
				\omega_n(f)=&-\frac{1}{2 \pi i} \oint_{\Gamma} \frac{(\upsilon_{X}-1)f(z)\underline{m}_n^0(z)^2} {z\left(1-{c_n} \int \underline{m}_n^0(z)^2 t^2(1+t \underline{m}_n^0(z))^{-2} d H_n(t)\right)}\\&\qquad\qquad\times\Bigg(\frac{{c_n} P_n(z)\int \underline{m}_n^0(z) t^2(1+t \underline{m}_n^0(z))^{-3} d H_n(t)}{\left(1-{c_n}\int \underline{m}_n^0(z)^2 t^2(1+t \underline{m}_n^0(z))^{-2} d H_n(t)\right)} -Q_n(z)\Bigg)dz \\
				&-\frac{1}{2 \pi i} \oint_{\Gamma} \mu_{X}f(z)z^2\underline{m}_n^0(z)^2\Bigg[\underline{m}_n^0(z)P_n(z)\widetilde{U}_n^1(z,z)\\&\qquad\qquad\times\bigg(1-{c_n}\int\frac{\underline{m}_n^0(z)^2t^2dH_n(t)}{(1+\underline{m}_n^0(z)t)^2}\bigg)^{-1}-\widetilde{V}_n^1(z,z)\Bigg]dz
			\end{split}
		\end{align}
		and $\bm{\Omega_1}$ is an $r\times r$ matrix with the $(s,t)$th entry being
		\begin{equation}\label{Varsample}
			\left(\bm{\Omega_1}\right)_{st}
			=-\frac{1}{4 \pi^2} \iint_{\Gamma_1\times \Gamma_2} f_{s}(z_1)f_{t}(z_2)\lim_{n\rightarrow\infty}\left( \upsilon_{X}C_n^1(z_1,z_2)+\mu_{X}C_n^2(z_1,z_2)\right)d z_1 d z_2.
		\end{equation}
		The functions $C_n^1$, $C_n^2$ are expressed as
		\begin{align}\label{Cov2sam}
			\begin{split}
				C_n^1(z_1,z_2)=&\frac{(\underline{m}_2-\underline{m}_1)z_1z_2}{z_2-z_1}\bigg(V_n^3(z_1,z_2)+z_2^2\underline{m}_2^2g_n(z_2)V_n^2(z_2,z_1)\\&+z_1^2\underline{m}_1^2g_n(z_1)V_n^2(z_1,z_2)+z_1^2z_2^2\underline{m}_1^2\underline{m}_2^2g_n(z_1)g_n(z_2)U_n^2(z_1,z_2)\bigg)\\&+\frac{(\underline{m}_2-\underline{m}_1)^2z_1z_2}{\underline{m}_1\underline{m}_2(z_2-z_1)^2}\bigg(z_1z_2\underline{m}_1\underline{m}_2\zeta_n^1(z_1,z_2)\zeta_n^1(z_2,z_1)\\&-z_1\underline{m}_1g_n(z_1)\zeta_n^1(z_1,z_2)-z_2\underline{m}_2g_n(z_2)\zeta_n^1(z_2,z_1)+g_n(z_1)g_n(z_2)a_n(z_1,z_2)\bigg),
			\end{split}
		\end{align}
		and
		\begin{align}\label{65sam}
			\begin{split}
				&C_n^2(z_1,z_2)=z_1z_2\underline{m}_1\underline{m}_2\bigg(\widetilde{V}_n^3(z_1,z_2)+z_2^2\underline{m}_2^2g_n(z_2)\widetilde{V}_n^2(z_2,z_1)+z_1^2\underline{m}_1^2g_n(z_2)\widetilde{V}_n^2(z_1,z_2)\\
                &+z_2^2\underline{m}_2^2g_n(z_2)z_1^2\underline{m}_1^2g_n(z_1)\widetilde{U}_n^2(z_1,z_2)-z_1\underline{m}_1g_n(z_1)\widetilde{V}_n^1(z_1,z_2)-z_1\underline{m}_1z_2^2\underline{m}_2^2g_n(z_1)g_n(z_2)\widetilde{U}_n^1(z_1,z_2)\\&-z_2\underline{m}_2g_n(z_2)\widetilde{V}_n^1(z_2,z_1)-z_2\underline{m}_2z_1^2\underline{m}_1^2g_n(z_1)g_n(z_2)\widetilde{U}_n^1(z_2,z_1)+g_n(z_1)g_n(z_2)\tilde{a}_n(z_1,z_2)
				\bigg).
			\end{split}
		\end{align}
{Here $\underline{m}_i$ denotes $\underline{m}(z_i)$ for simplicity and the other n-associated terms are defined in detail in (\ref{za}).} The contours $\Gamma_1$ and $\Gamma_2$ are disjoint and have the same properties as $\Gamma$.	
	\end{theorem}
	
We look at the special case when $\bm{B}_n=\bm{I}_n$. Obviously it satisfies Assumption \ref{asc} (i) since $k_n=n$ now. It can be easily checked that $\frac{1}{n}\operatorname{tr}(z\bm{I}_n+z\underline{m}_n^0(z)\bm{\Sigma}_n)^{-1}\bm{B}_n=\int\frac {dH_n(t)}{z(1+\underline{m}_n^0(z)t)}={m}_n^0(z)$. Then $\Theta_n(f)$ in equation \eqref{add2-1} reduces to
$$\Theta_n(f)=n\int f(x)d\left(F^{\bm{S}_n}(x)-F^{c_n,H_n}(x)\right),$$
 which is the conventional linear spectral statistic corresponding to the sample covariance matrix (see \citep{Bai:2004}). And our theoretical result in Theorem \ref{main} coincides with the traditional one (see our Remark C.1 for detailed calculations).

The asymptotic covariances (\ref{Varsample}) are mainly determined by two functions $C_n^1(z_1,z_2)$ and $C_n^2(z_1,z_2)$ defined in (\ref{Cov2sam}) and (\ref{65sam}).  If the first four moments of the underlying distribution matches with that of a standard Gaussian distribution, then $\mu_{X}=0$ and  $C_n^2(z_1,z_2)$ disappears in \eqref{Varsample}.  Our Remark C.1 shows certain cases that  the $n$-associated terms in (\ref{Cov2sam}) and (\ref{65sam}) are convergent and have succinct forms. Moreover, it can be seen from our proof that these terms are uniformly bounded in $z\in\mathcal{C}$ (see (C.3)), where $\mathcal{C}$ is any contour in the complex plane enclosing the closed interval (\ref{interval}). Therefore, in application, we often use a normalized version of Theorem \ref{main}, which is summarized in the following Proposition \ref{normalize}.

	\begin{proposition}\label{normalize}
		Suppose Assumptions \ref{asa}, \ref{asb} and \ref{asc} (i) hold. We further assume that $\lambda_{r}(\bm{\Omega}_n^1)\geq c_1>0$ for large $n$ and some positive constant $c_1$, where
		$$\left(\bm{\Omega_n^1}\right)_{st}
		=-\frac{1}{4 \pi^2} \iint_{\Gamma_1\times \Gamma_2} f_{s}(z_1)f_{t}(z_2)\left( \upsilon_{X}C_n^1(z_1,z_2)+\mu_{X}C_n^2(z_1,z_2)\right)d z_1 d z_2.$$
		Then we have
		\begin{align}\label{converge4}
			(\bm{\Omega_n^1})^{-1/2}(\Theta_n(f_1)-\omega_n(f_1),\cdots,\Theta_n(f_r)-\omega_n(f_r))^{\top}\stackrel{D}{\rightarrow}\mathcal{N}(\bm{0},\bm{I}_r).
		\end{align}
	\end{proposition}
    \begin{remark}\label{remark11}
	     The condition $\lambda_{r}(\bm{\Omega}_n^1)\geq c_1>0$ actually implies the linearly independence of $f_1,\cdots,f_r$ in the sense that for any unit vector $\bm{u}\in \mathbb{R}^r$, the variance of $\Theta_n\left((f_1,\cdots,f_r)\bm{u}\right)$ does not approach $0$.
      \end{remark}


The asymptotic distribution of GLSS is then investigated when $k_n=o(n)$. It should be noted that when $k_n/n\rightarrow 0$, the quantities relevant to $\bm{B}_n$ in Theorem \ref{main} all become zeros, resulting in $\Theta_n(f)\stackrel{\mathbb{P}}{\rightarrow} 0 $. Consequently, we need to seek for a suitable sequence $a_n \rightarrow\infty$, such that $a_n\Theta_n(f)$ converges to a non-degenerate distribution.
	\begin{theorem}\label{the6.1} [$k_n$ is much smaller than $n$].
Suppose that Assumptions \ref{asa}, \ref{asb} and \ref{asc} (ii) hold. Define
  $$H_n^1(z_1,z_2)=\frac{1}{k_n}\operatorname{tr}\left(\bm{B}_n\bm{\overline{\Sigma}}_n^{-1}(z_1)\bm{\Sigma}_n\bm{\overline{\Sigma}}_n^{-1}(z_2)\right)^2,$$
  and
  {\footnotesize  $$H_n^2(z_1,z_2)=\frac{{\underline{m}_n^0(z_1)\underline{m}_n^0(z_2)}}{k_nz_1z_2}\sum_{i=1}^n\left(\bm{\Sigma}_n^{1/2}\bm{\overline{\Sigma}}_n^{-1}(z_1)\bm{B}_n\bm{\overline{\Sigma}}_n^{-1}(z_1)\bm{\Sigma}_n^{1/2}\right)_{ii}\left(\bm{\Sigma}_n^{1/2}\bm{\overline{\Sigma}}_n^{-1}(z_2)\bm{B}_n\bm{\overline{\Sigma}}_n^{-1}(z_2)\bm{\Sigma}_n^{1/2}\right)_{ii}.$$
  }
		Then we have
		
		(i) the random vector
		\begin{equation}\label{Vector1}
			\sqrt{\frac{N}{k_n}}\left(\Theta_n(f_1), \ldots, \Theta_n(f_r)\right)
		\end{equation}
		forms a tight sequence in $n$.
		
		(ii) The random vector (\ref{Vector1}) converges weakly to a mean-zero $r$-dimensional Gaussian distribution, i.e.,
		\begin{align}\label{converge1}
			\begin{split}
				\sqrt{\frac N{k_n}}(\Theta_n(f_1),\cdots,\Theta_n(f_r))\stackrel{D}{\rightarrow}\mathcal{N}(\bm{0},\bm{\Omega_2}),
			\end{split}
		\end{align}	
		where $\bm{\Omega_2}$ is an $r\times r$ matrix with the $(s,t)$th entry being
{\small		\begin{align}\label{Var3}
			\begin{split}
				(\bm{\Omega_2})_{st}
				=-\frac{1}{4 \pi^2} \iint_{\Gamma_1\times \Gamma_2} f_{s}(z_1)f_{t}(z_2) \lim_{n\rightarrow\infty}\left(\frac{\upsilon_X(\underline{m}(z_2)-\underline{m}(z_1))H_n^1(z_1,z_2)}{z_1z_2(z_2-z_1)}+\mu_X H_n^2(z_1,z_2)\right)dz_1dz_2 ,
			\end{split}
		\end{align}
}
		where $\Gamma_1, \Gamma_2$ are assumed to be disjoint as described in Theorem \ref{main}.
	\end{theorem}

	\begin{remark}\label{xrmk2} Our Theorem \ref{the6.1} generalizes the results in \citep{Bai:2007} and \citep{Pan:2008}, which specifically consider the case $k_n=1$. To elaborate,
by assuming $\bm{B}_n=\bm{b}_n\bm{b}_n^*$, \citep{Bai:2007} obtained the CLT under Gaussian-like fourth moment assumption, i.e., $\mathbb{E}|X_{11}|^4=3$ in the real case and $\mathbb{E}|X_{11}|^4=2$ in the complex case. Moreover, they require
    \begin{align}\label{Bai}
        \sqrt{N}\left|\bm{b}_n^*\bm{\overline{\Sigma}}_n^{-1}(z)\bm{b}_n-\frac1n\operatorname{tr}\bm{\overline{\Sigma}}_n^{-1}(z)\right|\rightarrow0.
    \end{align}
 Condition \eqref{Bai} implies the convergence of $H_n^1(z_1,z_2)$, as confirmed by the equality (4.25) in \citep{Bai:2007}. \citep{Pan:2008} extended the moment condition to $\mathbb{E}|X_{11}|^4<\infty$, and additionally they required both \eqref{Bai} and
    \begin{align}\label{Pan}
        \max_{i}\left|\bm{b}_n^*\bm{\overline{\Sigma}}_n^{-1}(z_1)\bm{\Sigma}_n^{1/2}\bm{e}_i\right|\rightarrow 0.
    \end{align}
    It is evident that \eqref{Pan} directly indicates $H_n^2(z_1,z_2)\rightarrow 0$. Therefore within their specified frameworks, our result established in Theorem \ref{the6.1} aligns with theirs.
	\end{remark}
	
\begin{remark}\label{rmkadd} In application, we may also use a normalized version of Theorem \ref{the6.1} as done in Proposition \ref{normalize}. Comparing Theorem \ref{the6.1} with Theorem \ref{main}, one can see that the asymptotic mean in (\ref{converge1}) is zero, which is totally different from that in (\ref{converge}) where a bias  $\omega_n(f)$ appears. Also, the expression for asymptotic variance is significantly simplified when $k_n=o(n)$ compared to the case when $k_n/n\geq c_0$.
\end{remark}

\begin{remark}\label{remarkprok}
      {
        In Section G of the supplementary material, we  establish counterparts to Theorems \ref{main} and \ref{the6.1} using the L\'{e}vy-Prohorov distance (see \citep{Najim:2015}), thereby  removing the conditions $c_n \rightarrow c$ in Assumption \ref{asa} and $H_n \rightarrow H$ in Assumption \ref{asb}. }
    \end{remark}
    
We give a further illustration on the non-random part $\frac{1}{2\pi i}\oint_{\Gamma}f(z)\operatorname{tr}(z\bm{I}_n+z\underline{m}_n^0(z)\bm{\Sigma}_n)^{-1}\bm{B}_ndz$ in $\Theta_{n}(f)$.  Analogous to expression (\ref{sep}), it can be rewritten as
	\begin{equation}\label{add-nonrand}
 \frac{1}{2\pi i}\sum_{i=1}^{k_n}\sum_{j=1}^n|\langle\bm{b}_i,\bm{\upsilon}_j\rangle|^2\oint_{\Gamma}\frac{s_if(z)}{z\left(1+\lambda_{j}(\bm{\Sigma}_n)\underline{m}_n^0(z)\right)}dz,
 \end{equation}
	where the decomposition $\bm{\Sigma_n}=\sum_{j=1}^n\lambda_{j}(\bm{\Sigma_n})\bm{\upsilon}_j\bm{\upsilon}_j^*$ is employed. Each summation term in \eqref{add-nonrand} is divided into two parts: one determined by the inner product of the eigenvectors of $\bm{B}_n$ and $\bm{\Sigma}_n$, and the other solely influenced by the eigenvalues. Consequently, if the inner product $\langle\bm{b}_i,\bm{\upsilon}_j\rangle=0$ for some $i,j$, then the corresponding summation term becomes zero. The non-random part \eqref{add-nonrand} is governed by the non-orthogonal eigenvectors of $\bm{B}_n$ and $\bm{\Sigma}_n$. Therefore, it is possible for us to design a suitable GLSS for a specified hypothesis testing regarding the eigenspace structure, as exemplified in Section \ref{section5}.

{A careful examination of our  Theorems \ref{main} and \ref{the6.1} reveals that their statements can be unified into a single expression, which we summarize as follows.
\begin{theorem}\label{unify}
    Suppose Assumptions \ref{asa}, \ref{asb} hold and $\tau_{\bm{B}_n}=k_n/n\rightarrow\tau\in[0,1]$. Recall the definitions in (\ref{za}) and define $\mathcal{P}_n(z)=N/k_nP_n(z)$, $\mathcal{Q}_n(z)=N/k_nQ_n(z)$, $g_n^1(z)=N/k_ng_n(z)$, $\tilde{\zeta}_n^1(z_1,z_2)=N/k_n\zeta_n^1(z_,z_2)$, $\mathcal{V}_n^i(z_1,z_2)=N/k_nV_n^i(z_1,z_2)$, $\widetilde{\mathcal{V}}_n^i(z_1,z_2)=N/k_n$ $\widetilde{V}_n^i(z_1,z_2)$ for $i=1,2,3$. We have the following results:

    (i) the random vector
		\begin{equation}\label{newvector}
			\sqrt{N/k_n}\left(\Theta_n(f_1), \ldots, \Theta_n(f_r)\right)
		\end{equation}
		forms a tight sequence in $n$.
		
		(ii) After suitable centralization, the random vector (\ref{newvector}) converges weakly to an $r$-dimensional Gaussian distribution, i.e.,
		\begin{align}\label{convergenew}
			\sqrt{N/k_n}(\Theta_n(f_1)-\omega_n(f_1),\cdots,\Theta_n(f_r)-\omega_n(f_r))\stackrel{D}{\rightarrow}\mathcal{N}(\bm{0},\bm{\Omega}_3),
		\end{align}
		where $\omega_n(f)$ is defined in (\ref{meansam}) and $\bm{\Omega}_3$ is an $r\times r$ matrix with the $(s,t)$th entry being
		\begin{equation}\label{Varsamplenew}
			\left(\bm{\Omega}_3\right)_{st}
			=-\frac{1}{4 \pi^2} \iint_{\Gamma_1\times \Gamma_2} f_{s}(z_1)f_{t}(z_2)\lim_{n\rightarrow\infty}\left( \upsilon_{X}\mathcal{C}_n^1(z_1,z_2)+\mu_{X}\mathcal{C}_n^2(z_1,z_2)\right)d z_1 d z_2.
		\end{equation}
		The functions $\mathcal{C}_n^1$, $\mathcal{C}_n^2$ are expressed as
		\begin{align}\label{Cov2samnew}
			\begin{split}
				\mathcal{C}_n^1(z_1,z_2)=&\frac{(\underline{m}_2-\underline{m}_1)z_1z_2}{z_2-z_1}\bigg(\mathcal{V}_n^3(z_1,z_2)+c\tau z_2^2\underline{m}_2^2g_n^1(z_2)\mathcal{V}_n^2(z_2,z_1)\\&+c\tau z_1^2\underline{m}_1^2g_n^1(z_1)\mathcal{V}_n^2(z_1,z_2)+c\tau z_1^2z_2^2\underline{m}_1^2\underline{m}_2^2g_n^1(z_1)g_n^1(z_2)U_n^2(z_1,z_2)\bigg)\\&+\frac{c\tau(\underline{m}_2-\underline{m}_1)^2z_1z_2}{\underline{m}_1\underline{m}_2(z_2-z_1)^2}\bigg(z_1z_2\underline{m}_1\underline{m}_2\tilde{\zeta}_n^1(z_1,z_2)\tilde{\zeta}_n^1(z_2,z_1)\\&-z_1\underline{m}_1g_n^1(z_1)\tilde{\zeta}_n^1(z_1,z_2)-z_2\underline{m}_2g_n^1(z_2)\tilde{\zeta}_n^1(z_2,z_1)+g_n^1(z_1)g_n^1(z_2)a_n(z_1,z_2)\bigg),
			\end{split}
		\end{align}
		and
		\begin{align}\label{65samnew}
			\begin{split}
				\frac{\mathcal{C}_n^2(z_1,z_2)}{z_1z_2\underline{m}_1\underline{m}_2}=&\widetilde{\mathcal{V}}_n^3(z_1,z_2)+c\tau z_2^2\underline{m}_2^2g_n^1(z_2)\widetilde{\mathcal{V}}_n^2(z_2,z_1)+c\tau z_1^2\underline{m}_1^2g_n^1(z_2)\widetilde{\mathcal{V}}_n^2(z_1,z_2)\\
                &+c\tau z_2^2\underline{m}_2^2g_n^1(z_2)z_1^2\underline{m}_1^2g_n^1(z_1)\widetilde{U}_n^2(z_1,z_2)-c\tau z_1\underline{m}_1g_n^1(z_1)\widetilde{\mathcal{V}}_n^1(z_1,z_2)\\&-c\tau z_2\underline{m}_2g_n^1(z_2)\widetilde{\mathcal{V}}_n^1(z_2,z_1)-c\tau z_2\underline{m}_2z_1^2\underline{m}_1^2g_n^1(z_1)g_n^1(z_2)\widetilde{U}_n^1(z_2,z_1)\\&-c\tau z_1\underline{m}_1z_2^2\underline{m}_2^2g_n^1(z_1)g_n^1(z_2)\widetilde{U}_n^1(z_1,z_2)+c\tau g_n^1(z_1)g_n^1(z_2)\tilde{a}_n(z_1,z_2).
			\end{split}
		\end{align}
    Here $\underline{m}_i$ denotes $\underline{m}(z_i)$ for simplicity. The contours $\Gamma_1$ and $\Gamma_2$ are disjoint, as described in Theorem \ref{main}.	
\end{theorem}}

\subsection{Two explicit examples}	\label{secexam}
{In this section, we present two explicit examples -- one for the special case
$\bm{\Sigma}_n = \bm{I}$ and another for a general $\bm{\Sigma}_n$ -- and derive closed-form expressions for the corresponding limiting distributions in both cases.  In each example, the matrix $\bm{B}_n$ is taken to be general, subject only to Assumption \ref{asb}.
    }
    \begin{example}\label{Example1}
      {Consider a special $\bm{\Sigma}_n=\bm{I}_n$.  Suppose that Assumptions \ref{asa}, \ref{asb} hold, and that $\tau_{\bm{B}_n}=k_n/n\rightarrow\tau\in[0,1]$. Let $f_k(x)=x^k, k=1,2$, and $f_3(x)=\log(x)$. We further denote $\mu_{\bm{B}_n}=\operatorname{tr}(\bm{B}_n)/k_n$, $s_{\bm{B}_n}=\operatorname{tr}(\bm{B}_n^2)/k_n$ and $d_{\bm{B}_n}=\sum_{i=1}^n(\bm{B}_n)_{ii}^2/k_n$.
\begin{itemize}
\item  If $c\in(0,1)$, we can establish the joint distribution of $\operatorname{tr}f_1(\bm{S}_n)\bm{B}_n$, $\operatorname{tr}f_2(\bm{S}_n)\bm{B}_n$, and $\operatorname{tr}f_3(\bm{S}_n)\bm{B}_n$ as follows:
    \begin{align*}
        \sqrt{n/k_n}\begin{pmatrix} \operatorname{tr}f_1(\bm{S}_n)\bm{B}_n-(nq_1+p_1)\tau_{\bm{B}_n}\mu_{\bm{B}_n} \\\operatorname{tr}f_2(\bm{S}_n)\bm{B}_n-(nq_2+p_2)\tau_{\bm{B}_n}\mu_{\bm{B}_n}\\\operatorname{tr}f_3(\bm{S}_n)\bm{B}_n-(nq_3+p_3)\tau_{\bm{B}_n}\mu_{\bm{B}_n}
\end{pmatrix}\stackrel{D}{\rightarrow}\mathcal{N}(\bm{0},\lim_{n\rightarrow\infty}\bm{\Omega}_n^{(3)}),
    \end{align*}
    where
    \begin{align*}
        q_1=1, \quad q_2=c_n+1, \quad q_3=\frac{c_n-1}{c_n}\log(1-c_n)-1,
    \end{align*}
    \begin{align*}
    \begin{split}
        &p_1=0,\quad p_2=c(\upsilon_X-1+\mu_X),\quad p_3=\frac{(\upsilon_X-1)\log(1-c)}{2}-\frac{c\mu_X}{2},
    \end{split}
    \end{align*}
    and $\bm{\Omega}_{n}^{(3)}$ is a $3\times 3$ matrix with
    \small{\begin{align*}
        (\bm{\Omega}_{n}^{(3)})_{i,j}=\upsilon_X r_1(f_i,f_j)s_{\bm{B}_n}+\upsilon_X r_2(f_i,f_j)\mu_{\bm{B}_n}^2\tau+\mu_X r_3(f_i,f_j)d_{\bm{B}_n}+\mu_X r_4(f_i,f_j)\mu_{\bm{B}_n}^2\tau.
    \end{align*}}
    Here, in the definition of $(\bm{\Omega}_{n}^{(3)})_{i,j}$, the explicit expressions for each $r_k(f_i,f_j)$ ($k=1,2,3,4$) are given as follows:
    \small{\begin{align*}
    \begin{split}
        &r_1(f_1,f_1)=c, \quad r_1(f_1,f_2)=c^2+2c,\quad r_1(f_2,f_2)=c^3+5c^2+4c,\quad \\& r_1(f_1,f_3)=\frac c2+\frac{1-c}c\log(1-c)+1,\quad r_1(f_2,f_3)=\frac {c^2}6+\frac{5c}2+\frac{1-c^2}c\log(1-c)+1,\\&r_1(f_3,f_3)=1+2\operatorname{Li}_2(c)+\frac{c^2-1}{c^2}[\log(1-c)]^2, \ \text{with}\ \operatorname{Li}_2(c)=\sum_{j=1}^{\infty}c^j/j^2,
        \end{split}
    \end{align*}}
     \small{\begin{align*}
    \begin{split}
        &r_2(f_1,f_1)=0, \quad r_2(f_1,f_2)=c^2, \quad r_2(f_2,f_2)=3c^3+5c^2,\\&r_2(f_1,f_3)=\frac c2-\frac{1-c}c\log(1-c)-1,\quad r_2(f_2,f_3)=\frac {5c^2}6-\frac{c}2-\frac{1-c^2}c\log(1-c)-1,\\&r_2(f_3,f_3)=-\log(1-c)-1-2\operatorname{Li}_2(c)-\frac{c^2-1}{c^2}[\log(1-c)]^2, \ \text{with}\ \operatorname{Li}_2(c)=\sum_{j=1}^{\infty}c^j/j^2,
        \end{split}
    \end{align*}}
    \small{\begin{align*}
    \begin{split}
        &r_3(f_1,f_1)=c, \quad r_3(f_1,f_2)=c^2+2c, \quad r_3(f_2,f_2)=c^3+4c^2+4c,\\& r_3(f_1,f_3)=\frac c2+\frac{1-c}c\log(1-c)+1,\quad r_3(f_2,f_3)=\frac {c^2}2+2c+\frac{(1-c)(c+2)}c\log(1-c)+2,\\& r_3(f_3,f_3)=c\left(\frac{c-1}{c^2}\log(1-c)-\frac12-\frac1c\right)^2,
        \end{split}
    \end{align*}}
    and
    \small {\begin{align*}
    \begin{split}
        &r_4(f_1,f_1)=0, \quad r_4(f_1,f_2)=c^2, \quad r_4(f_2,f_2)=3c^3+4c^2,\\& r_4(f_1,f_3)=\frac c2-\frac{1-c}c\log(1-c)-1,\quad r_4(f_2,f_3)=\frac {3c^2}2-\frac{(1-c)(c+2)}c\log(1-c)-2,\\& r_4(f_3,f_3)=c-c\left(\frac{c-1}{c^2}\log(1-c)-\frac12-\frac1c\right)^2.
        \end{split}
    \end{align*}}
\item In general, when $c\in(0,\infty)$, $f_3(x)=\log(x)$ may be undefined. In this case, we establish the joint distribution of $\operatorname{tr}f_1(\bm{S}_n)\bm{B}_n$ and $\operatorname{tr}f_2(\bm{S}_n)\bm{B}_n)^T$ as follows:
     \begin{align*}
        \sqrt{n/k_n}\begin{pmatrix} \operatorname{tr}f_1(\bm{S}_n)\bm{B}_n-(nq_1+p_1)\tau_{\bm{B}_n}\mu_{\bm{B}_n} \\\operatorname{tr}f_2(\bm{S}_n)\bm{B}_n-(nq_2+p_2)\tau_{\bm{B}_n}\mu_{\bm{B}_n}
\end{pmatrix}\stackrel{D}{\rightarrow}\mathcal{N}(\bm{0},\lim_{n\rightarrow\infty}\bm{\Omega}_n^{(3)}\big|_{1:2}^{1:2}),
    \end{align*}
    where $\bm{\Omega}_n^{(3)}\big|_{1:2}^{1:2}$ is the $2\times2$ upper-left sub-matrix of $\bm{\Omega}_n^{(3)}$.
\end{itemize}}
    \end{example}
    {The proof of Example \ref{Example1} is deferred to Section H.1 of the supplementary material. Actually, in Section H.1, we consider the first four power functions $x,x^2,x^3,x^4$, and $\log(x)$ (see Theorem H.1). However, the asymptotic variances of $\operatorname{tr}\bm{S}_n^3\bm{B}_n$ and $\operatorname{tr}\bm{S}_n^4\bm{B}_n$, as well as their covariances with $\operatorname{tr}\log(\bm{S}_n)\bm{B}_n$, are quite lengthy. Therefore, due to space constraints, we do not include the cases $x^3$ and $x^4$ in the main paper.}

    \begin{example}\label{Example2}
   Consider a general $\bm{\Sigma}_n$.
Suppose that Assumptions 2.1 and 2.2 hold, and that $\tau_{\bm{B}_n}=k_n/n\rightarrow\tau\in[0,1]$.    Define $\bm{\Sigma}_B^{k,s}=\bm{\Sigma}^{k/2}\bm{B}\bm{\Sigma}^{s/2}$, where, for simplicity, we use $\bm{B}$ and $\bm{\Sigma}$ to denote $\bm{B}_n$ and $\bm{\Sigma}_n$, respectively.  We have
    \begin{align}\label{Aos50}
        \sqrt{N/k_n}\binom{\operatorname{tr} \bm{S}_n \bm{B} -\operatorname{tr}\bm{\Sigma}\bm{B}}{\operatorname{tr} \bm{S}_n^2 \bm{B}-\mu_{n,2}}\xrightarrow{D} \mathcal{N}\left(\bm{0},\lim_{n\rightarrow\infty} \begin{pmatrix}   \upsilon_X\sigma_{11,n}^{(1)}+\mu_X\sigma_{11,n}^{(2)}&\upsilon_X\sigma_{12,n}^{(1)}+\mu_X\sigma_{12,n}^{(2)}\\\upsilon_X\sigma_{12,n}^{(1)}+\mu_X\sigma_{12,n}^{(2)}&\upsilon_X\sigma_{22,n}^{(1)}+\mu_X\sigma_{22,n}^{(2)}
\end{pmatrix}\right),
    \end{align}
    where
    \begin{align}\label{Aos51}
        \mu_{n,2}=\left(1+N^{-1}(\upsilon_X-1)\right)\operatorname{tr}\bm{\Sigma}^2\bm{B}+N^{-1}\operatorname{tr}\bm{\Sigma}\operatorname{tr}\bm{\Sigma}\bm{B}+\mu_XN^{-1}\operatorname{tr}\left(\bm{\Sigma}\circ\bm{{\Sigma}}_{B}^{1,1}\right),
    \end{align}
    \begin{align*}
        \sigma_{11,n}^{(1)}=k_n^{-1}\operatorname{tr}\bm{\Sigma}\bm{B}\bm{\Sigma}\bm{B},\qquad\sigma_{11,n}^{(2)}=k_n^{-1}\operatorname{tr}\left(\bm{\Sigma}_B^{1,1}\circ\bm{\Sigma}_B^{1,1}\right),
    \end{align*}
    \begin{align*}
        \sigma_{12,n}^{(1)}=2k_n^{-1}\operatorname{tr}\bm{\Sigma}^2\bm{B}\bm{\Sigma}\bm{B}+k_n^{-1}N^{-1}\operatorname{tr}\bm{\Sigma}\bm{B}\bm{\Sigma}\bm{B}\cdot\operatorname{tr}\bm{\Sigma}+c\tau k_n^{-2}\operatorname{tr}\bm{\Sigma}^2\bm{B}\cdot\operatorname{tr}\bm{\Sigma}\bm{B},
    \end{align*}
    \begin{align*}
        \begin{split}\sigma_{12,n}^{(2)}=&k_n^{-1}N^{-1}\operatorname{tr}\left(\bm{\Sigma}_B^{1,1}\circ\bm{\Sigma}_B^{1,1}\right)\operatorname{tr}\bm{\Sigma}+c\tau k_n^{-2}\operatorname{tr}\left(\bm{\Sigma}\circ\bm{\Sigma}_B^{1,1}\right)\operatorname{tr}\bm{\Sigma}\bm{B}\\&+k_n^{-1}\operatorname{tr}\left(\bm{\Sigma}_B^{1,3}\circ\bm{\Sigma}_B^{1,1}\right)+k_n^{-1}\operatorname{tr}\left(\bm{\Sigma}_B^{3,1}\circ\bm{\Sigma}_B^{1,1}\right),
        \end{split}
    \end{align*}
    \begin{align}\label{Aos55}
        \begin{split}\sigma_{22,n}^{(1)}=&2k_n^{-1}\operatorname{tr}\bm{\Sigma}^3\bm{B}\bm{\Sigma}\bm{B}+2k_n^{-1}\operatorname{tr}\bm{\Sigma}^2\bm{B}\bm{\Sigma}^2\bm{B}+4k_n^{-1}N^{-1}\operatorname{tr}\bm{\Sigma}^2\bm{B}\bm{\Sigma}\bm{B}\cdot\operatorname{tr}\bm{\Sigma}\\&+4c\tau k_n^{-2}\operatorname{tr}\bm{\Sigma}^3\bm{B}\cdot\operatorname{tr}\bm{\Sigma}\bm{B}+c\tau k_n^{-2}\left(\operatorname{tr}\bm{\Sigma}^2\bm{B}\right)^2+k_n^{-1}N^{-1}\operatorname{tr}\bm{\Sigma}\bm{B}\bm{\Sigma}\bm{B}\cdot\operatorname{tr}\bm{\Sigma}^2\\&+k_n^{-1}N^{-2}\operatorname{tr}\bm{\Sigma}\bm{B}\bm{\Sigma}\bm{B}\cdot\left(\operatorname{tr}\bm{\Sigma}\right)^2+c\tau k_n^{-2}N^{-1}\left(\operatorname{tr}\bm{\Sigma}\bm{B}\right)^2\operatorname{tr}\bm{\Sigma}^2\\&+2c\tau k_n^{-2}N^{-1}\operatorname{tr}\bm{\Sigma}\bm{B}\cdot\operatorname{tr}\bm{\Sigma}^2\bm{B}\cdot\operatorname{tr}\bm{\Sigma},
        \end{split}
    \end{align}
    and
    \begin{align}\label{Aos56}    \begin{split}\sigma_{22,n}^{(2)}=&k_n^{-1}\operatorname{tr}\left(\bm{\Sigma}_B^{1,3}\circ\bm{\Sigma}_B^{1,3}\right)+k_n^{-1}\operatorname{tr}\left(\bm{\Sigma}_B^{3,1}\circ\bm{\Sigma}_B^{3,1}\right)+2k_n^{-1}\operatorname{tr}\left(\bm{\Sigma}_B^{1,3}\circ\bm{\Sigma}_B^{3,1}\right)\\&+2k_n^{-1}N^{-1}\operatorname{tr}\bm{\Sigma}\left[\operatorname{tr}\left(\bm{\Sigma}_B^{3,1}\circ\bm{\Sigma}_B^{1,1}\right)+\operatorname{tr}\left(\bm{\Sigma}_B^{1,3}\circ\bm{\Sigma}_B^{1,1}\right)\right]\\&+2c\tau k_n^{-2}\operatorname{tr}\bm{\Sigma}\bm{B}\cdot\operatorname{tr}\left(\bm{\Sigma}_B^{3,1}\circ\bm{\Sigma}\right)+2c\tau k_n^{-2}\operatorname{tr}\bm{\Sigma}\bm{B}\cdot\operatorname{tr}\left(\bm{\Sigma}_B^{1,3}\circ\bm{\Sigma}\right)\\&+k_n^{-1}N^{-2}\left(\operatorname{tr}\bm{\Sigma}\right)^2\operatorname{tr}\left(\bm{\Sigma}_B^{1,1}\circ\bm{\Sigma}_B^{1,1}\right)+c\tau k_n^{-2}N^{-1}\left(\operatorname{tr}\bm{\Sigma}\bm{B}\right)^2\operatorname{tr}\left(\bm{\Sigma}\circ\bm{\Sigma}\right)\\&+2c\tau k_n^{-2}N^{-1}\operatorname{tr}\bm{\Sigma}\cdot\operatorname{tr}\bm{\Sigma}\bm{B}\cdot\operatorname{tr}\left(\bm{\Sigma}_B^{1,1}\circ\bm{\Sigma}\right).
 \end{split}
    \end{align}
    \end{example}

{The proof of Example \ref{Example2} is provided in Section H.2 of the supplementary material.}

	\section{Simulations}\label{section4}
	\color{black}In this section, a series of simulations are conducted with varying choices of $\bm{\Sigma}_n$, $\bm{B}_n$ and underlying distributions of $X_{ij}$ to empirically validate the theoretical results presented in Section \ref{section2}. In Section \ref{caseB1}, we choose $\operatorname{rank}(\bm{B}_n)=n$ to satisfy the conditions stated in Theorem \ref{main}, while in Section \ref{caseB2} we select some constant values for $\operatorname{rank}(\bm{B}_n)$ that align with Theorem \ref{the6.1}. Let $r=1$, $f(z)=z^2$ and {$n=500$, $N=1000$. We consider the real case, which means $\upsilon_X=2$}. Denote
	$$\widetilde{\Theta}_n(f)=\left(2\sigma_{22,n}^{(1)}+\mu_X\sigma_{22,n}^{(2)}\right)^{-1/2}\sqrt{N/k_n}\left(\operatorname{tr}\bm{S}_n^2\bm{B}_n-\mu_{n,2}\right),$$
    where $\mu_{n,2}$, $\sigma_{22,n}^{(1)}$, and $\sigma_{22,n}^{(2)}$ are defined in (\ref{Aos51}), (\ref{Aos55}), and (\ref{Aos56}) respectively. \color{black}
 Our theoretical findings suggest that the distribution of $\widetilde{\Theta}_n(f)$ converges to $\mathcal{N}(0,1)$. All numerical results presented below are based on {$M=5000$} replications, yielding {5000} simulated estimates $\left(\widetilde{\Theta}_n^1(f),\cdots,\widetilde{\Theta}_n^M(f)\right)$ of $\widetilde{\Theta}_n(f)$. The empirical mean and variance are
	\begin{equation}\label{120}
		\widehat{\mathbb{E}X_f}=\frac1M \sum_{k=1}^M\widetilde{\Theta}_n^k(f),
	\end{equation}
	and
	\begin{equation}\label{121}
		\widehat{\operatorname{Var}X_f}=\frac1M \sum_{k=1}^M(\widetilde{\Theta}_n^k(f)-\widehat{\mathbb{E}X_f})^2.
	\end{equation}
  Besides, we compute the following quantities
\begin{align}\label{simulatedalpha}
	\hat{\alpha}_r=\frac1M\sum_{k=1}^M\mathbb{I}_{\{\widetilde{\Theta}_n^k(f)>\Phi^{-1}(1-\alpha)\}},\qquad \hat{\alpha}_l=\frac1M\sum_{k=1}^M\mathbb{I}_{\{\widetilde{\Theta}_n^k(f)<\Phi^{-1}(\alpha)\}},
\end{align}
where $\Phi$ is the distribution function of $\mathcal{N}(0,1)$. In the following simulations we fix $\alpha=0.05$.

\color{black}

Eight different models will be considered. For each model, we plot the histogram of $\left(\widetilde{\Theta}_n^1(f),\cdots,\widetilde{\Theta}_n^M(f)\right)$ and compare it with the density function of $\mathcal{N}(0,1)$. Additionally, the normal QQ-plot is presented to further validate the asymptotical normality.

	\subsection{{The matrix \texorpdfstring{$\textbf{B}_n$}{} is of full rank}}\label{caseB1}
This section will investigate six distinct models, each offering different choices of $\bm{\Sigma}_n$ and $\bm{B}_n$, as well as varying underlying distributions of $X_{ij}$. In all these models, $\bm{B}_n$ possesses full rank, which aligns with the condition stated in Theorem \ref{main}.

\textbf{Model 1. $\bm{\Sigma}_n=\bm{I}_n$}, $X_{ij}\sim \mathcal{N}(0,1)$ and $\bm{B}_n$ is a diagonal matrix with the $i$-th entry being $(i/n+1)$.

	\textbf{Model 2. $\bm{\Sigma}_n=\bm{I}_n$}, $X_{ij}\sim (\text{Gamma}(2,1)-2)/\sqrt{2}$ and $\bm{B}_n$ is a diagonal matrix with the $i$-th entry being $(i/n+1)$. Model 2 differs from Model 1 in the way of selecting the distribution of $X_{ij}$. In Model 2, $X_{ij}$ follows a gamma distribution with shape parameter being $2$ and scale parameter being $1$. We subtract $2$ and divide by $\sqrt{2}$ to ensure $\mathbb{E}|X_{11}|^2=1$. One can easily check that $\mathbb{E}|X_{11}|^4=6$, different from that of $\mathcal{N}(0,1)$.
	
\color{black}	\textbf{Model 3.} $\bm{\Sigma}_n$ is the covariance matrix of $AR(1)$ sequence with coefficient $0.5$ (i.e. the $(i,j)$th entry is $0.5^{|i-j|}$), $X_{ij}\sim \mathcal{N}(0,1)$ and $\bm{B}_n=\bm{\Sigma}_n$.
	
	\textbf{Model 4.} $\bm{\Sigma}_n$ is a diagonal matrix with $(\bm{\Sigma}_n)_{ii}=(i/n)^2+0.2$, {$X_{ij}\sim (\text{Gamma}(2,1)-2)/\sqrt{2}$} and $\bm{B}_n$ is a diagonal matrix with $(\bm{B}_n)_{ii}=i/n+0.2$.
	
	\textbf{Model 5.} $\bm{\Sigma}_n$ is the same as in Model 3, $X_{ij}\sim \mathcal{N}(0,1)$ and $\bm{B}_n$ is chosen to be an arbitrary realization of the standard Wigner matrix.

	\textbf{Model 6.} $\bm{\Sigma}_n$ and $\bm{B}_n$ are the same as in Model 5, and {$X_{ij}\sim (\text{Gamma}(2,1)-2)/\sqrt{2}$} whose fourth moment is different from that of $\mathcal{N}(0,1)$.
\color{black}	

{The histogram plots and QQ plots are depicted in Figures \ref{tt1}-\ref{tt2} for Models 1-2 and in Figures K.1-K.4 in Section K of the supplementary material for Models 3-6, respectively}. These results confirm the accuracy of our theoretical results.

	\begin{figure}[H]		
		\centering
		\begin{subfigure}{0.48\linewidth}
			\centering
			\caption{Histogram}
			\includegraphics[width=0.8\linewidth]{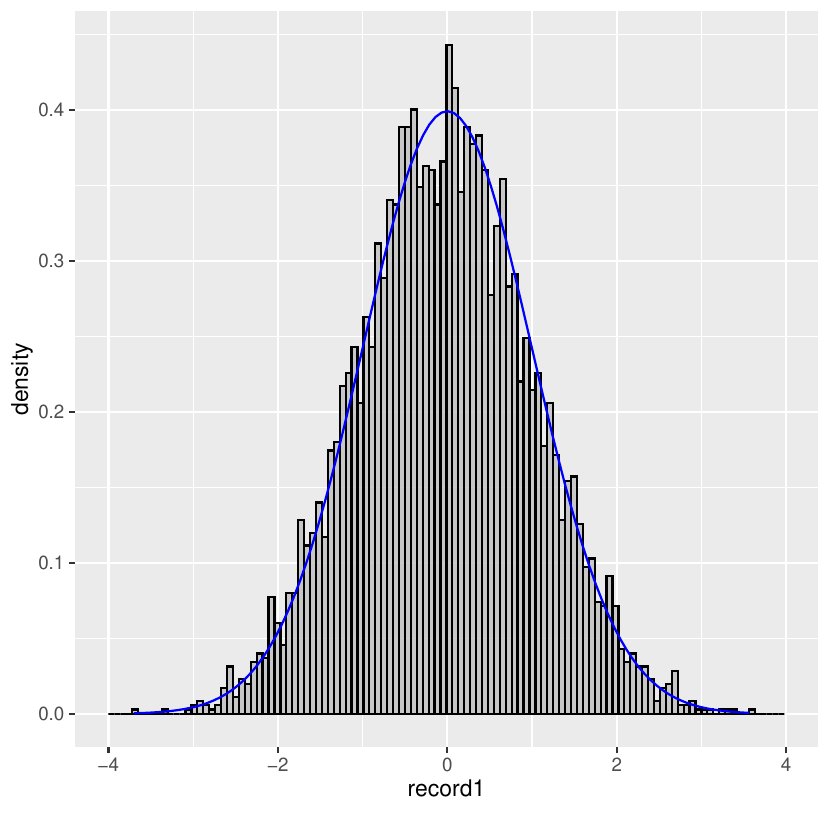}
			\label{t1}
		\end{subfigure}
		\begin{subfigure}{0.48\linewidth}
			\centering
			\caption{QQ-plot }
			\includegraphics[width=0.8\linewidth]{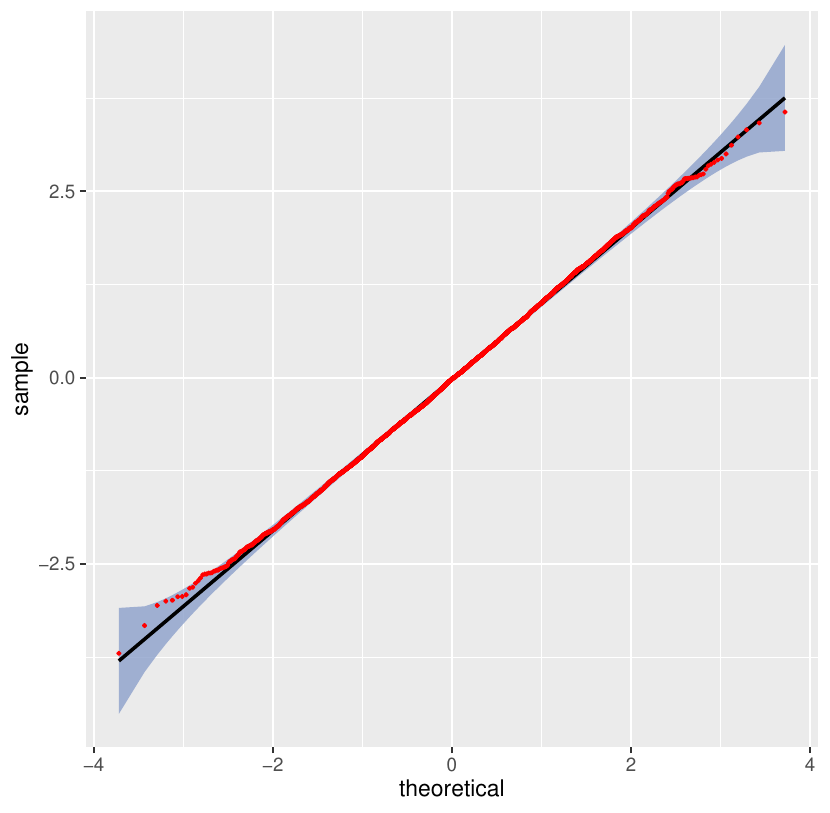}
			\label{t2}
		\end{subfigure}
		\centering
		\caption{Model 1: (a): Histogram of the records $\left(\widetilde{\Theta}_n^1(f),\cdots,\widetilde{\Theta}_n^M(f)\right)$ {with $X_{ij}\sim\mathcal{N}(0,1)$} and density curve of $\mathcal{N}(0,1)$ (blue line) (b):   QQ-plot of the records.}
		\label{tt1}
	\end{figure}

	\begin{figure}[H]		
		\centering
		\begin{subfigure}{0.48\linewidth}
			\centering
			\caption{Histogram}
			\includegraphics[width=0.8\linewidth]{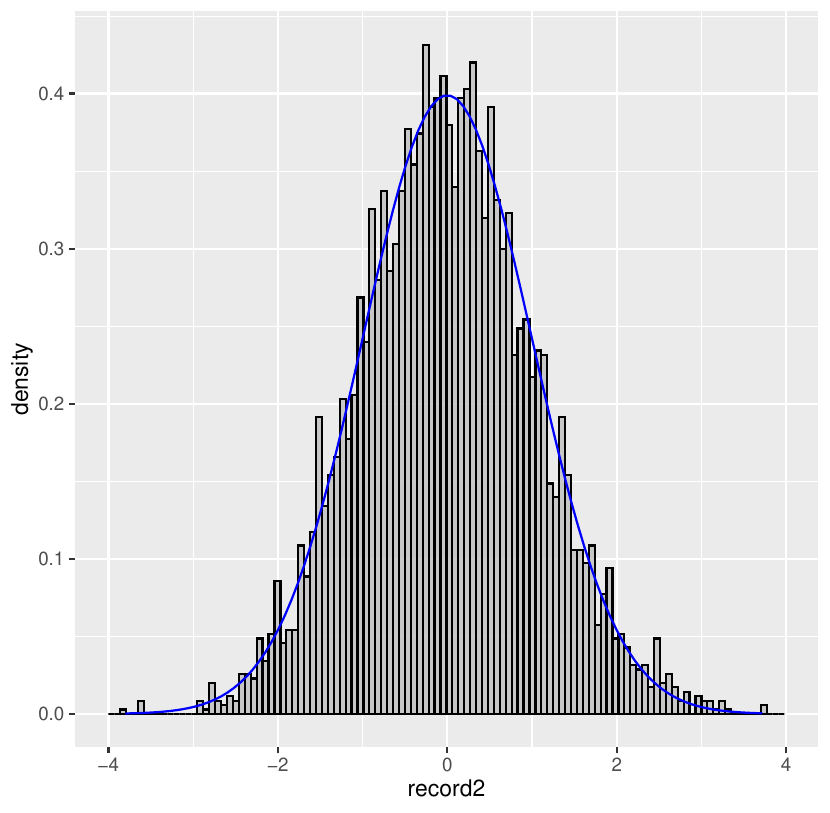}
			\label{t3}
		\end{subfigure}
		\begin{subfigure}{0.48\linewidth}
			\centering
			\caption{QQ-plot }
			\includegraphics[width=0.8\linewidth]{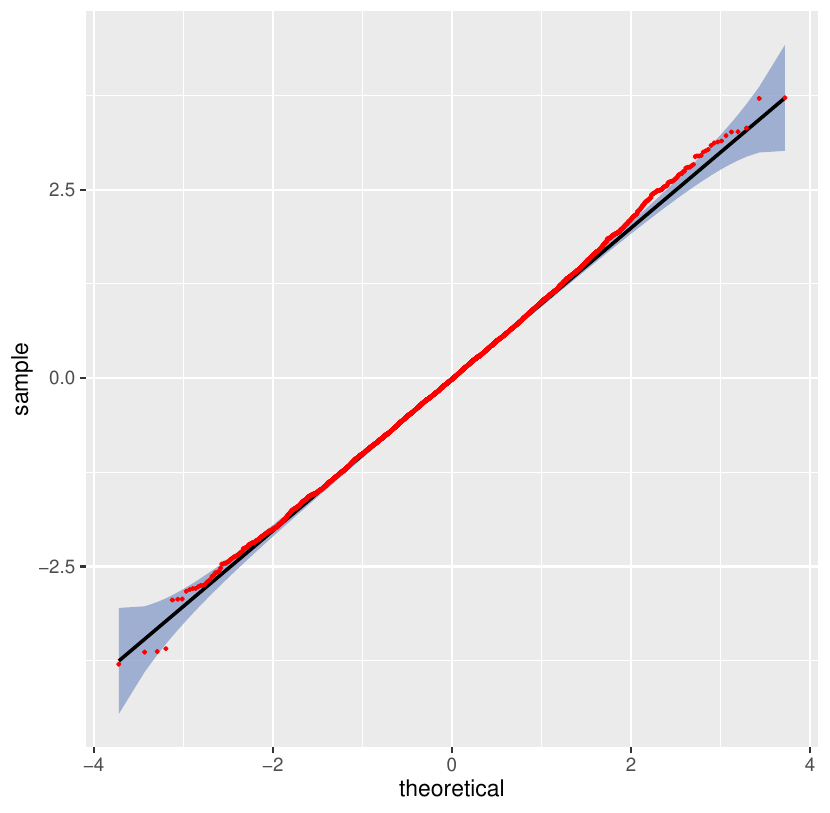}
			\label{t4}
		\end{subfigure}
		\centering
		\caption{Model 2: (a): Histogram of the  records $\left(\widetilde{\Theta}_n^1(f),\cdots,\widetilde{\Theta}_n^M(f)\right)$ {with $X_{ij}\sim (\text{Gamma}(2,1)-2)/\sqrt{2}$} and density curve of $\mathcal{N}(0,1)$ (blue line) (b):   QQ-plot of the records.}
		\label{tt2}
	\end{figure}

	\subsection{{The matrix \texorpdfstring{$\textbf{B}_n$}{} is of low rank}}\label{caseB2}
For the models considered in this section, the ranks of $\bm{B}_n$ are constant values, which aligns with the condition stated in Theorem \ref{the6.1}.
	
	\textbf{Model 7.} $\bm{\Sigma}_n=\bm{I}_n$, $X_{ij}\sim \mathcal{N}(0,1)$ and $\operatorname{rank}(\bm{B}_n)=5$. Specifically, $\bm{B}_n$ is a diagonal matrix with $(\bm{B}_n)_{ii}=i/2$, for $i=1,\cdots,5$.
	
	\textbf{Model 8.} $\bm{\Sigma}_n$ is the same as in Model 3, {$X_{ij}\sim (\text{Gamma}(2,1)-2)/\sqrt{2}$} and $\operatorname{rank}(\bm{B}_n)=10$. Specifically, $\bm{B}_n=\sum_{i=1}^{10}\bm{b}_i\bm{b}_i^*$, where $\bm{b}_i$'s are selected from the eigenvectors of a realization for Wigner matrix.
	
\begin{figure}[H]		
		\centering
		\begin{subfigure}{0.48\linewidth}
			\centering
			\caption{Histogram}
			\includegraphics[width=0.8\linewidth]{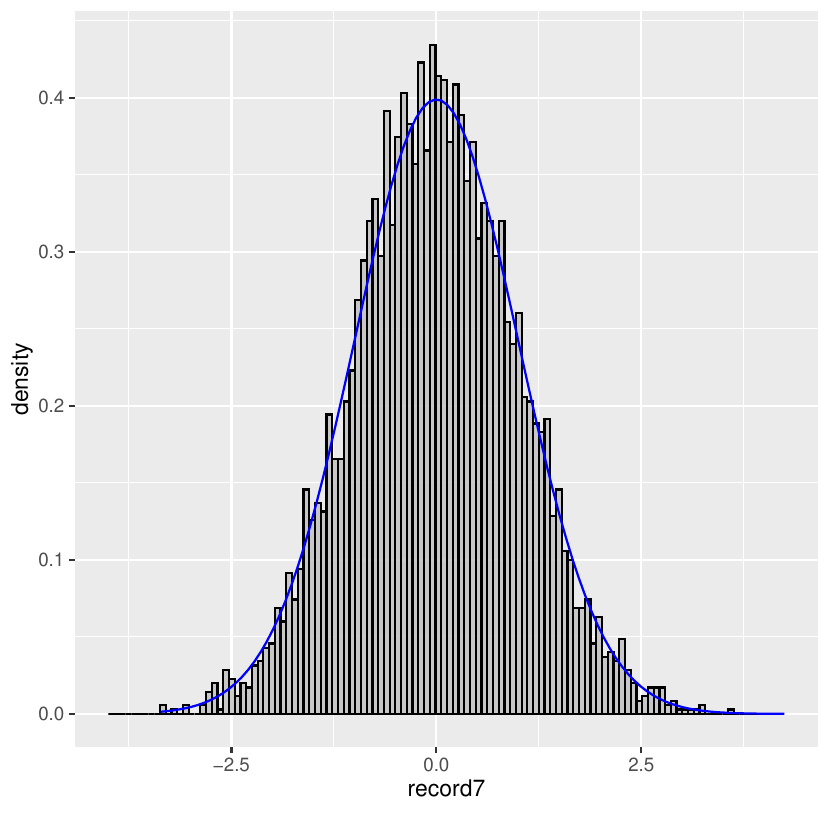}
			\label{t13}
		\end{subfigure}
		\begin{subfigure}{0.48\linewidth}
			\centering
			\caption{QQ-plot }
			\includegraphics[width=0.8\linewidth]{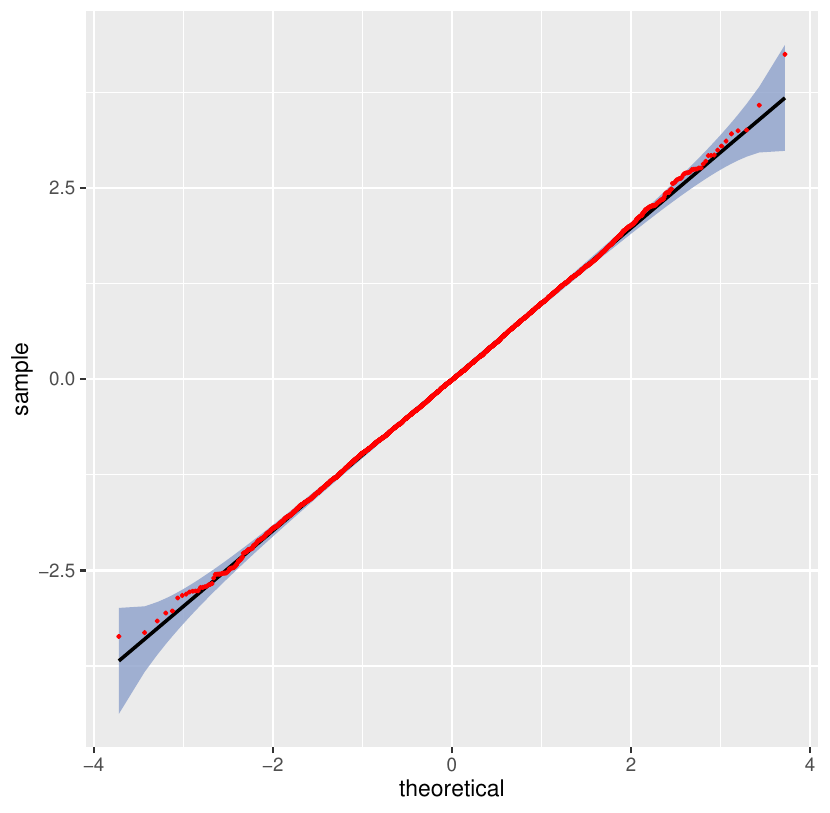}
			\label{t14}
		\end{subfigure}
		\centering
		\caption{Model 7: (a): Histogram of the records $\left(\widetilde{\Theta}_n^1(f),\cdots,\widetilde{\Theta}_n^M(f)\right)^{\top}$ {with $X_{ij}\sim\mathcal{N}(0,1)$} and density curve of $\mathcal{N}(0,1)$ (blue line) (b):   QQ-plot of the records.}
		\label{tt7}
	 \end{figure}

\begin{figure}[H]		
		\centering
		\begin{subfigure}{0.48\linewidth}
			\centering
			\caption{Histogram}
			\includegraphics[width=0.9\linewidth]{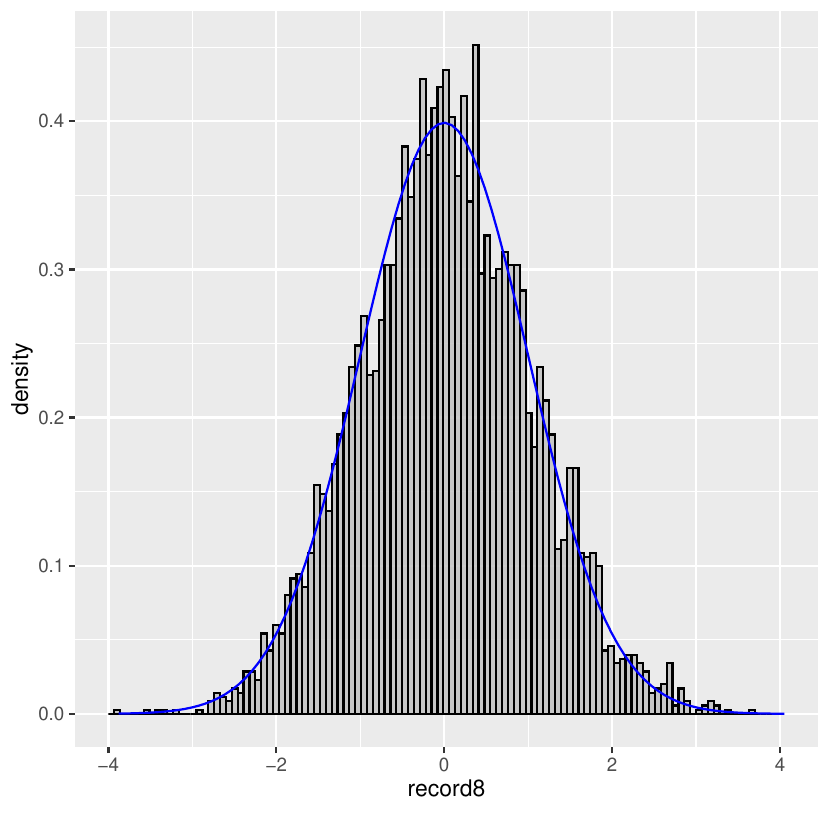}
			\label{t17}
		\end{subfigure}
		\begin{subfigure}{0.48\linewidth}
			\centering
			\caption{QQ-plot }
			\includegraphics[width=0.9\linewidth]{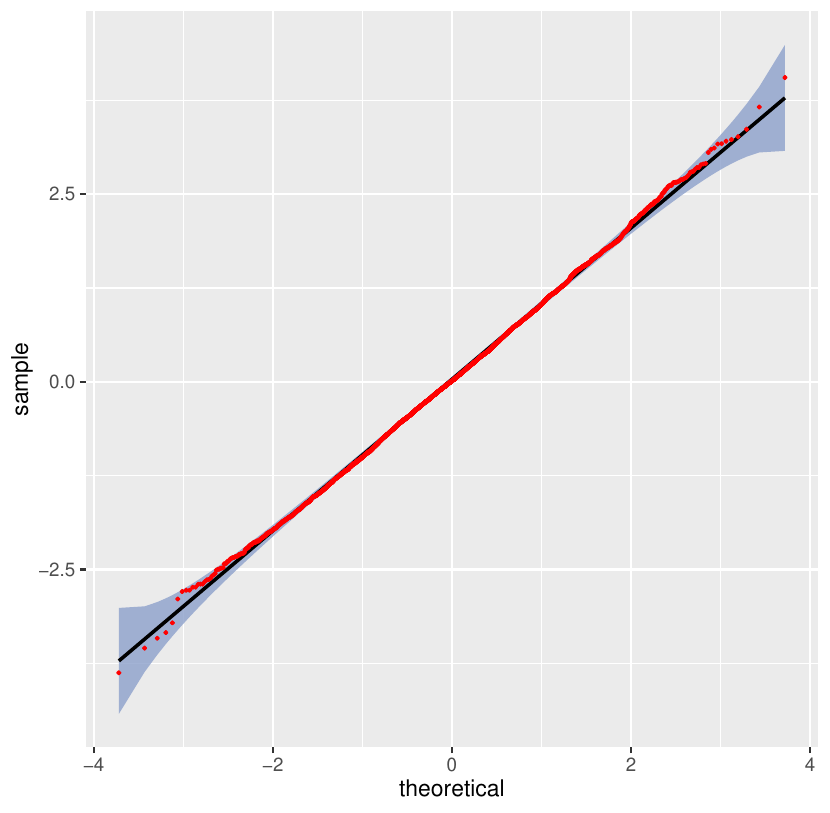}
			\label{t18}
		\end{subfigure}
		\centering
		\caption{Model 8: (a): Histogram of the records $\left(\widetilde{\Theta}_n^1(f),\cdots,\widetilde{\Theta}_n^M(f)\right)^{\top}$ {with $X_{ij}\sim (\text{Gamma}(2,1)-2)/\sqrt{2}$} and density curve of $\mathcal{N}(0,1)$ (blue line) (b):   QQ-plot of the records.}
		\label{tt9}
	\end{figure}

Figures \ref{tt7} and \ref{tt9} present the histograms and QQ plots for the above two models, which demonstrate the accuracy of our theoretical results in Theorem \ref{the6.1}.
	
Finally, the empirical means and variances calculated by \eqref{120} and \eqref{121} for the aforementioned eight models are recorded in Table \ref{Table1}. {Besides, the quantities $\hat{\alpha}_r$ and $\hat{\alpha}_l$ defined in \eqref{simulatedalpha} are also reported in Table \ref{Table1}.} It is evident that, across all models, the empirical mean closely approximates zero while the variance closely approximates one {and the quantities are close to $0.05$}, thereby providing strong support for our theoretical findings.

	\begin{table}[H]
		\begin{center}
			\caption{Empirical mean and variance defined in (\ref{120}) and (\ref{121}) and the quantities defined in (\ref{simulatedalpha}) for the eight different models.}
			\label{Table1}
			\renewcommand\arraystretch{2}
			\begin{tabular}{|c|c|c|c|c|c|c|c|c|}
				\hline   \textbf{Model} & \textbf{1} & \textbf{2} & \textbf{3} & \textbf{4} &\textbf{5}&\textbf{6}&\textbf{7}&\textbf{8}\\
				\hline  $\widehat{\mathbb{E}X_f}$ & \color{black}-0.0213 &\color{black} 0.0008 & \color{black}-0.0229 &\color{black} -0.0168 &\color{black}0.0015&\color{black}-0.0050&\color{black}-0.0017&\color{black}0.0247\\
				\hline   $\widehat{\operatorname{Var}X_f}$ & \color{black}1.0362&\color{black} 1.0349 &\color{black} 0.9774 &\color{black} 1.0361 &\color{black}0.9912 &\color{black}1.0019&\color{black}0.9798&\color{black}1.0299\\
				\hline\color{black}$\hat{\alpha}_r$ & \color{black}0.0528 &\color{black} 0.0562 & \color{black}0.0466 &\color{black} 0.0510 &\color{black}0.0516&\color{black}0.0522&\color{black}0.0474&\color{black}0.0570\\
				\hline \color{black}  $\hat{\alpha}_l$ &\color{black} 0.0560&\color{black} 0.0472 &\color{black} 0.0516 &\color{black} 0.0534 &\color{black}0.0498 &\color{black}0.0490&\color{black}0.0460&\color{black}0.0488\\
				\hline
			\end{tabular}
		\end{center} 	
	\end{table}
\color{black}

\section{Application to Eigenspace Testing on ``Population-spiked'' Covariance Matrices}\label{section5}

Hypothesis testing for eigenspaces of the spiked covariance matrix plays a crucial role in statistical machine learning and is encountered in various modern algorithms, see \citep{Silin:2020} for an extensive discussion on this topic. However, many existing methods for such problems are limited to the case when $n\ll N$, both theoretically and practically, unless there are constraints on the structure of the covariance matrix. See for example, bootstrap based approach \citep{Naumov:2019,Silin:2020}, Bayesian or Frequentist-Bayes related method \citep{Silin:2020,Silin:2018}, sample splitting method\citep{Koltchinskii:2017}, and the Le Cam optimal test proposed in \citep{Hallin:2010}. In the high-dimensional setting where $n\asymp N$, for the spiked covariance matrix model $\bm{\Sigma}_n$ that admits the decomposition
	\begin{equation}\label{300}
		\bm{\Sigma}_n=\bm{I}_n+\sum_{i=1}^{r_n}d_i\bm{v}_i\bm{v}_i^{\top}, \quad d_1\geq \cdots \geq d_{r_n}>0,
	\end{equation}
 \citep{Bao:2022} proposed a statistic based on the accurate results on the joint distribution of the few leading extreme eigenvalues and the generalized components of their associated eigenvectors. We would like to mention two assumptions required in \citep{Bao:2022}. Firstly, $r_n=r$ is a fixed constant. Secondly, their Assumption 2.4 imposes a restriction on the minimal distance of $|d_i-d_j|$ when $d_i\neq d_j$ and requires a positive lower bound  $\sqrt{c_n}$ for the spikes $d_i$, $i=1,\cdots, r$.

 In this section, we propose a novel approach based on GLSS to investigate the eigenspaces of covariance matrices exhibiting ``population-spiked'' characteristics. The term ``population-spiked'' is employed to distinguish our method from existing approaches in that it accommodates diverging number of spikes, while only requiring $0<\inf_n\min_{i=1,\cdots,r_n}d_i\leq\sup_n\max_{i=1,\cdots,r_n}d_i<\infty$ without imposing an additional positive lower bound for the magnitude of $d_i$.

%

\subsection{{Methodology and theoretical results}}
 We now present our methodology for testing whether the eigenspace spanned by the eigenvectors corresponding to the $r_n$ spikes is equivalent to a given subspace. Denote $\bm{\mathcal{Z}}_n=\sum_{i=1}^{r_n}\bm{v}_i\bm{v}_i^{\top}$. Then the testing problem is
	\begin{equation}\label{301}
		\boldsymbol{H}_0 :\quad \bm{\mathcal{Z}}_n=\bm{\mathcal{Z}}_0 \quad \operatorname{vs} \quad \boldsymbol{H}_1 :\quad \bm{\mathcal{Z}}_n\neq \bm{\mathcal{Z}}_0,
	\end{equation}
	for a given projection matrix $\bm{\mathcal{Z}}_0$. In the ideal case when $r_n/N\rightarrow 0$ and accurate estimation of all $d_i$'s at a rate of $\mathrm{o}_{\mathbb{P}}(N^{-1/2})$ is possible, Theorem \ref{the6.1} suggests a natural test statistic $\Theta_n(f)$ defined in \eqref{add2-1} by using $\bm{B}_n=\bm{\mathcal{Z}}_0$  for testing hypothesis \eqref{301}. However, it is practically impossible to achieve such an ideal estimator for $d_i$. Even when $r_n$ is fixed, according to Theorem 2.10 in \citep{Bao:2022}, the estimation of spiked eigenvalues exhibits robustness only up to a rate $\mathrm{O}_{\mathbb{P}}(N^{-1/2})$, not to mention when $r_n$ diverges.
In order to eliminate the effect of unknown ${d_i}$'s, we select $\bm{B}_n$ as the projection matrix orthogonal to $\bm{\mathcal{Z}}_0$, i.e. $\bm{B}_n=\bm{I}_n-\bm{\mathcal{Z}}_0$. Consequently, the $\operatorname{rank}(\bm{B}_n)$ now satisfies Assumption \ref{asc}(i) and Theorem \ref{main} implies a limiting Gaussian distribution for the test statistic $\Theta_n(f)$. Encouragingly, through this selection of $\bm{B}_n$, under the null hypothesis, neither the non-random component nor its asymptotic mean and variance in $\Theta_n(f)$ incorporate any unknown spiked eigenvalues. The sole remaining unknown term is $\underline{m}_n^0(z)$. Simply substituting $\underline{m}_n(z)$ for $\underline{m}_n^0(z)$ would impact the asymptotic distribution stated in Theorem \ref{main} due to an $\mathrm{O}_{\mathbb{P}}(N^{-1})$ order discrepancy between $\underline{m}_n^0(z)$ and $\underline{m}_n(z)$, which constitutes a non-negligible error. To surmount this challenge, we adapt $\Theta_n(f)$ by defining our test statistic as follows:
	\begin{equation}\label{302}
		\Delta_n(f)=\operatorname{tr}f(\bm{S}_n)(\bm{I}_n-\bm{\mathcal{Z}}_0)-\frac{n-r_n}{2\pi i}\oint_{\Gamma}\frac{f(z)}{z+z\underline{m}_n(z)}dz,
	\end{equation}
and refer to this testing approach as \textbf{{Functional Projection}}.
Focusing on the case of real variables, which is commonly encountered in practical applications, we establish the asymptotic distribution of $\Delta_n(f)$ as presented in the following Theorem \ref{theapp}.
	\begin{theorem}\label{theapp}
		Suppose that the population covariance matrix $\bm{\Sigma}_n$ admits the decomposition (\ref{300}). In addition to Assumption \ref{asa}, we further assume that
		$$(FP).\quad 0<\inf_n\min_{i=1,\cdots,r_n}d_i\leq\sup_n\max_{i=1,\cdots,r_n}d_i<\infty, \quad \operatorname{and} \quad r_n/N\rightarrow 0.$$
		Then under the null hypothesis $\boldsymbol{H}_0$ in \eqref{301}, we have
		\begin{equation}\label{304}
			\frac{\Delta_n(f)-\mu(f,r_n,n,N)}{\sqrt{\varrho(f,r_n,n,N)}}\stackrel{D}{\rightarrow} \mathcal{N}(0,1),
		\end{equation}
		where $\mu(f,r_n,n,N)$ and $\varrho(f,r_n,n,N)$ are explicitly defined by means of equations (1.1)-(1.19) in the supplementary material.
	\end{theorem}

 \begin{remark}\label{rmk42}
{\textbf{[Universality on $d_j$'s]}}.  Suppose that $\bm{S}_n$ owns the spectral decomposition $\bm{S}_n=\sum_{j=1}^n\lambda_j\bm{u}_j\bm{u}_j^{*}$.
It has been observed that $\bm{u}_j$ exhibits distinct asymptotic behaviors under two scenarios, namely $d_j>\sqrt{c_n}$ and $d_j<\sqrt{c_n}$ (refer to \citep{Blo:2016} for example). Theorem \ref{theapp} reveals an intriguing phenomenon that our statistic $\Delta_n(f)$ consistently follows an asymptotically normal distribution as long as $\min_{j=1,\cdots,r_n}d_j>0$.  This implies that, contrary to conventional wisdom, hypothesis test (\ref{301}) can still be conducted even when all $d_j$'s are less than $\sqrt{c_n}$. To empirically validate this finding, we perform empirical studies in Section \ref{addsimu} to check the efficiency of our functional projection approach (\ref{302}) when all $d_i$'s are smaller than $\sqrt{c_n}$. The simulated results displayed in Figure \ref{ttnew} support the universality of our functional projection approach against variations in $d_j$'s.


\end{remark}

  \begin{remark}\label{remark11a}
In practice, we need to estimate $d_i$ and $\underline{m}_n^0(z)$ in $ \mu(f,r_n,n,N)$ and $\varrho(f,r_n,n,N)$. A good estimator for $\underline{m}_n^0(z)$ is $\underline{m}_n(z)$ since $\underline{m}_n(z)-\underline{m}_n^0(z)=O_{\mathbb{P}}(N^{-1})$.
    Regarding $d_i$, we use a shrinkage estimator $\hat{d}_i$ to replace ${d}_i$:
    \begin{equation}\label{esti}
\hat{d}_i=\left\{
\begin{aligned}
&\frac12(-c_n-1+\lambda_i(\bm{S}_n))+\frac12\sqrt{(-c_n-1+\lambda_i(\bm{S}_n))^2-4c_n},\quad\lambda_i(\bm{S}_n)\geq (1+\sqrt{c_n})^2+\delta,\\
&0\qquad\qquad\qquad\qquad\qquad\qquad\qquad\qquad\qquad\qquad\qquad\quad,\quad\text{otherwise},\\
\end{aligned}
\right.
\end{equation}
where $\delta>0$ is any pre-specified constant. When $d_i>\sqrt{c_n}$ and is bounded away from infinity, it is verified from \citep{Bao:2022} that $\hat{d}_i$ is a consistent estimator for $d_i$ given a fixed $r_n$. Define $ \hat{\mu}(f,r_n,n,N)$ and $\hat{\varrho}(f,r_n,n,N)$ with $d_i$ and $\underline{m}_n^0(z)$ replaced by $\hat{d}_i$ and $\underline{m}_n(z)$. Since the $d_i$ associated terms in $\hat{\mu}(f,r_n,n,N)$ and $\hat{\varrho}(f,r_n,n,N)$ are of an order $O(r_n/N)$ (see eg. (1.1)), it is evident that $\hat{\varrho}(f,r_n,n,N)={\varrho}(f,r_n,n,N)+O_{\mathbb{P}}(r_n/N))$ and $\hat{\mu}(f,r_n,n,N)={\mu}(f,r_n,n,N)+O_{\mathbb{P}}(r_n/N))$ by the Assumption $(FP)$ in Theorem \ref{theapp}. As a consequence, (\ref{304}) still holds when ${\mu}(f,r_n,n,N)$ and ${\varrho}(f,r_n,n,N)$ are estimated by $\hat{\mu}(f,r_n,n,N)$ and $\hat{\varrho}(f,r_n,n,N)$.
  \end{remark}

	\subsection{{Numerical studies}}\label{sub5.1}
In this section, we conduct Monte Carlo simulations to investigate the finite-sample accuracy and power performance of our proposed testing approach--Functional Projection $\Delta_n(f)$ (abbreviated as \textbf{FP$\underline{~}\bm{f(z)}$}), and compare it with methods introduced in two existing papers. One is \citep{Silin:2020}, which utilized the bootstrapping method (abbreviated as \textbf{En$\underline{~}$Bo}) and the frequentist Bayes method (abbreviated \textbf{En$\underline{~}$Ba}) employing a power-enhanced norm with $s_1 =s_2 = 1$ (refer to their Definition 3.1).  We will use 1000 repetitions
	for both bootstrapping and frequentist Bayes procedure. The other one is the Fr-Adaptive (abbreviated as \textbf{Fr$\underline{~}$Ad}) proposed by \citep{Bao:2022}.

Without loss of generality, we assume that the eigenvectors align with the axes of the coordinate system under the null hypothesis $\boldsymbol{H}_0$, i.e. $\bm{v}_i=\bm{e}_i$ for $i=1,\cdots,r_n$. Then the hypothetical projection matrix is
	\begin{align}\label{HPM}
	\bm{\mathcal{Z}}_n=\left[\begin{array}{cc}
		\mathbf{I}_{r_n} & \mathbf{O}_{r_n\times(n-r_n)} \\
		\mathbf{O}_{(n-r_n) \times r_n} & \mathbf{O}_{(n-r_n) \times(n-r_n)}
	\end{array}\right]
	\end{align}
	and the default covariance matrix is diagonal with descending entries:
	\begin{align}\label{DCM}
	\bm{\Sigma}_n^{(0)}=\left[\begin{array}{llllll}
		1+d_1 & & & & & \\
		& \ddots & & & & \\
		& & 1+d_{r_n} & & & \\
		& & & 1 & & \\
		& & & & \ddots & \\
		& & & & & 1
	\end{array}\right].
	\end{align}
To study the performance under the alternative hypothesis, we follow the construction strategy in \citep{Silin:2020} and rotate the plane containing the first and the $(r_n+1)$-th axes by the angle $\varphi$, i.e. the leading eigenvector becomes
	$$
	\boldsymbol{v}_1^{\varphi}=[\underbrace{\cos \varphi, 0, \ldots, 0}_{r_n}, \sin \varphi, 0, \ldots, 0]^{\top} \text {, }
	$$
	while the $(r_n+1)$-th eigenvector turns into
	$$
	\boldsymbol{u}_{\varphi}=[\underbrace{-\sin \varphi, 0, \ldots, 0}_{r_n}, \cos \varphi, 0, \ldots, 0]^{\top} .
	$$
	The covariance matrix under $\boldsymbol{H}_1$ can be explicitly written as
	$$
	\bm{\Sigma}_n^{(\varphi)}=\left[\begin{array}{cccccccc}
		(1+d_1) \cos ^2 \varphi+ \sin ^2 \varphi & 0 & \ldots & 0 & d_1 \cos \varphi \sin \varphi & 0 & \ldots & 0 \\
		0 & 1+d_2 & & & 0 & & & \\
		\vdots & & \ddots & & \vdots & & & \\
		0 & & & 1+d_{r_n} & 0 & & & \\
		d_1\cos \varphi \sin \varphi & 0 & \ldots & 0 & (1+d_1)\sin ^2 \varphi+ \cos ^2 \varphi & & & \\
		0 & & & & & 1 & \\
		\vdots & & & & & & \ddots\\
		0 & & & & & & & 1
	\end{array}\right] .
	$$
A smaller $\varphi$ indicates a comparatively weaker alternative. The following scenarios will be taken into consideration.

	\textbf{$\bullet$ Scenario \uppercase\expandafter{\romannumeral1}.}
	Set $r_n=3$ with $d_1=9$, $d_2=5$ and $d_3=2$ (the spiked eigenvalues are simple with no multiplicity). The angle $\varphi$ varies within $\{1\%, 2\%, \cdots, 80\%\}\times \pi/2$ to capture the power performance trend. Both $X_{ij}\sim \mathcal{N}(0,1)$ and $X_{ij}\sim t(10)/\sqrt{5/4}$ are taken into account.
	
	\textbf{$\bullet$ Scenario \uppercase\expandafter{\romannumeral2}.}
	Set $d_1=9$ and $d_2=\cdots=d_{r_n}=4$ (eigenvalue multiplicity exists).  $X_{ij}\sim \mathcal{N}(0,1)$. Larger ranks $r_n=7$ and $r_n=11$ are considered. The angle $\varphi$ varies within $\{1\%, 2\%, \cdots, 80\%\}\times \pi/2$  to obtain the power performance trend.
	
		\textbf{$\bullet$ Scenario \uppercase\expandafter{\romannumeral3}.}
    Set $d_1=9$ and $d_2=\cdots=d_{r_n}=4$. $X_{ij}\sim \mathcal{N}(0,1)$. Fix $\varphi=\pi/8$ or $\varphi=0$, where the former reflects
    $\boldsymbol{H}_1$ and the latter corresponds to $\boldsymbol{H}_0$. The rank $r_n$ varies within $\{1,2,\cdots,15\}$ to check the tendency.

The choices for the remaining parameters are as follows: the nominal level $\alpha=0.1$, the threshold $\delta$ in (\ref{esti}) is $\delta=0.1$, the dimension $n=500$, the sample size $N\in\{500, 1000\}$, and the function $f(z)=z^2$ or $z^3$. The comparison of empirical powers is conducted using 100 replications, while the empirical sizes are calculated based on 1000 replications.

By setting $\varphi=0$, we record the empirical sizes in Scenarios I and II, as presented in Table \ref{table2}. It is observed that both our statistics FP$\underline{~}{z^2}$ and FP$\underline{~}{z^3}$ exhibit satisfactory accuracy, with the empirical size closely aligning with the nominal level $0.1$. In Scenario II, Fr$\underline{~}$Ad shows significantly inflated sizes, particularly when the number of spikes is large ($r_n=11$). Both En$\underline{~}$Bo and En$\underline{~}$Ba suffer from severe size distortion across all settings in Scenarios I and II.

\begin{figure}[H]		
		\centering
		\begin{subfigure}{0.48\linewidth}
			\centering
			\caption{n=500,N=500}
			\includegraphics[width=0.96\linewidth]{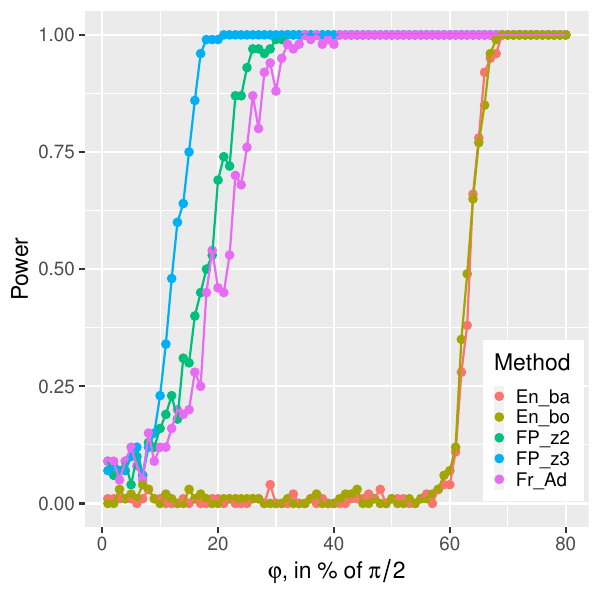}
			\label{t20}
		\end{subfigure}
		\begin{subfigure}{0.48\linewidth}
			\centering
			\caption{n=500,N=1000}
			\includegraphics[width=0.96\linewidth]{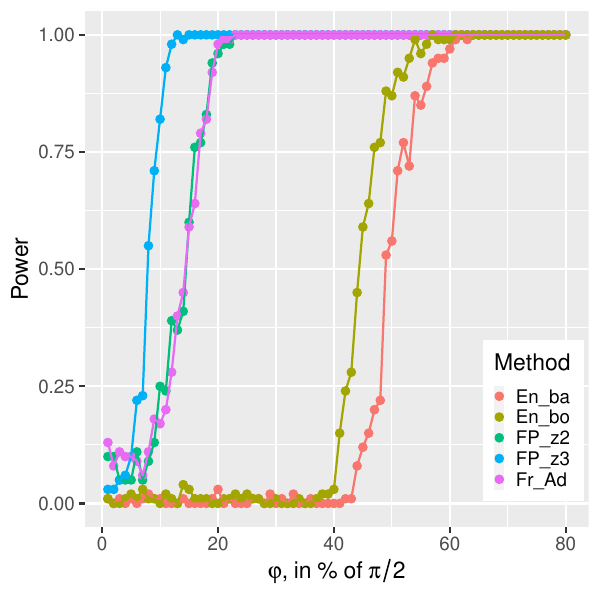}
			\label{t21}
		\end{subfigure}
		\centering
		\caption{Power comparison for Scenario I when $X_{ij}\sim \mathcal{N}(0,1)$. The angle $\varphi$ varies within $\{1\%, 2\%, \cdots, 80\%\}\times \pi/2$. The data dimension $n=500$. The sample size in the left plot (a) is $N=500$, while in the right plot (b) it is $N=1000$. FP$\underline{~}$z2 and FP$\underline{~}$z3 represents our approach FP with $f(z)=z^2$ and $z^3$, respectively.}
		\label{tt10}
	\end{figure}

\begin{figure}[H]		
		\centering
		\begin{subfigure}{0.48\linewidth}
			\centering
			\caption{n=500,N=500}
			\includegraphics[width=0.96\linewidth]{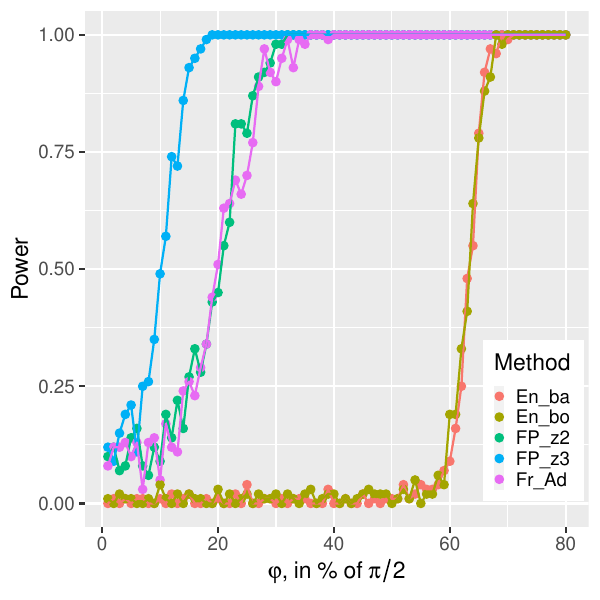}
			\label{t22}
		\end{subfigure}
		\begin{subfigure}{0.48\linewidth}
			\centering
			\caption{n=500,N=1000}
			\includegraphics[width=0.96\linewidth]{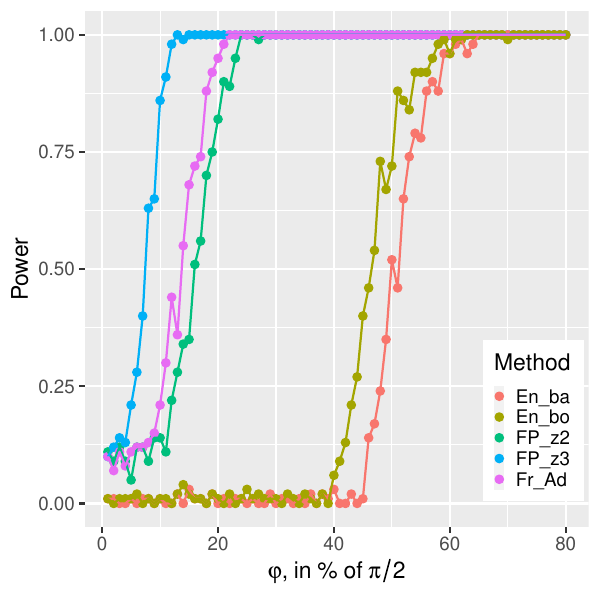}
			\label{t23}
		\end{subfigure}
		\centering
		\caption{Power comparison for Scenario I when $X_{ij}\sim t(10)/\sqrt{5/4}$. Others parameters are the same as introduced in Figure \ref{tt10}.}
		\label{tt11}
	\end{figure}

Figures \ref{tt10} and \ref{tt11} present the power comparison in Scenario I when $X_{ij}\sim \mathcal{N}(0,1)$ and $X_{ij}\sim t(10)/\sqrt{5/4}$, respectively. We can observe that our FP with $f(z)=z^3$ exhibits greater sensitivity and statistical power compared to other methods, particularly when the angle $\varphi$ is not large. The power of FP with $f(z)=z^2$ is comparable to that of Fr$\underline{~}$Ad. Both En$\underline{~}$Bo and En$\underline{~}$Ba show significantly reduced sensitivity to $\varphi$. This is evident from {the} observation that their power approaches 1 only within the range of $(\pi/2\times0.6,\pi/2\times0.8)$ for $\varphi$, while for smaller values of $\varphi$ than $\pi/2\times0.4$, the power remains close to zero.

\begin{figure}[H]		
		\centering
		\begin{subfigure}{0.48\linewidth}
			\centering
			\caption{n=500,N=500}
			\includegraphics[width=0.96\linewidth]{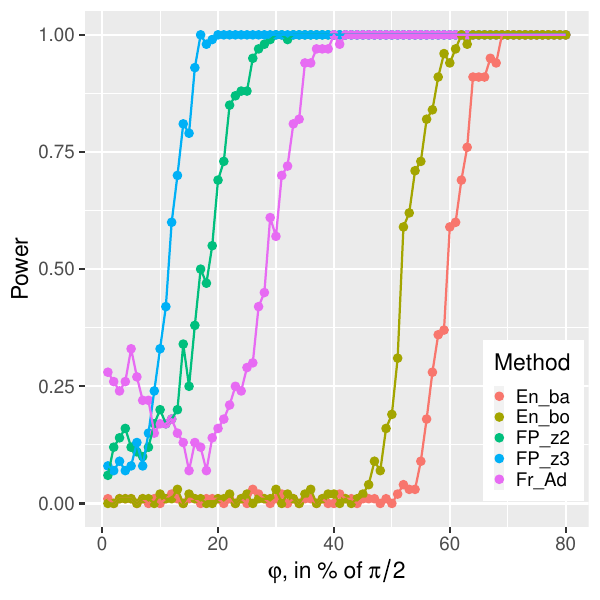}
			\label{t24}
		\end{subfigure}
		\begin{subfigure}{0.48\linewidth}
			\centering
			\caption{n=500,N=1000}
			\includegraphics[width=0.96\linewidth]{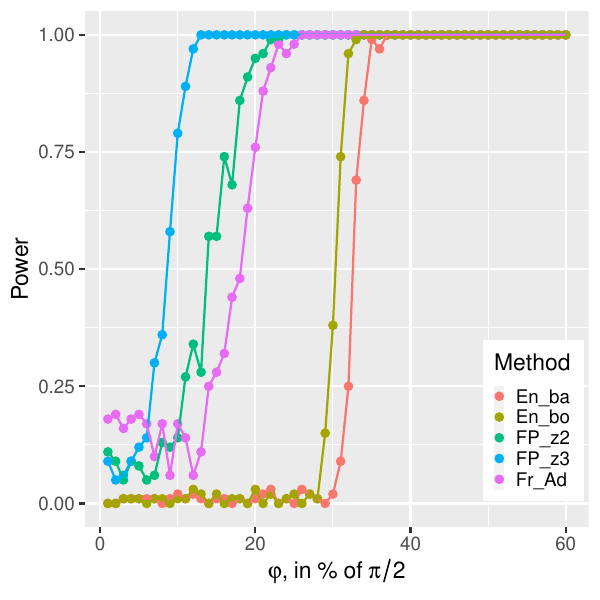}
			\label{t25}
		\end{subfigure}
		\centering
		\caption{Power comparison for Scenario II when $r_n=7$.}
		\label{tt12}
	\end{figure}

\begin{figure}[H]		
		\centering
		\begin{subfigure}{0.48\linewidth}
			\centering
			\caption{n=500,N=500}
			\includegraphics[width=0.96\linewidth]{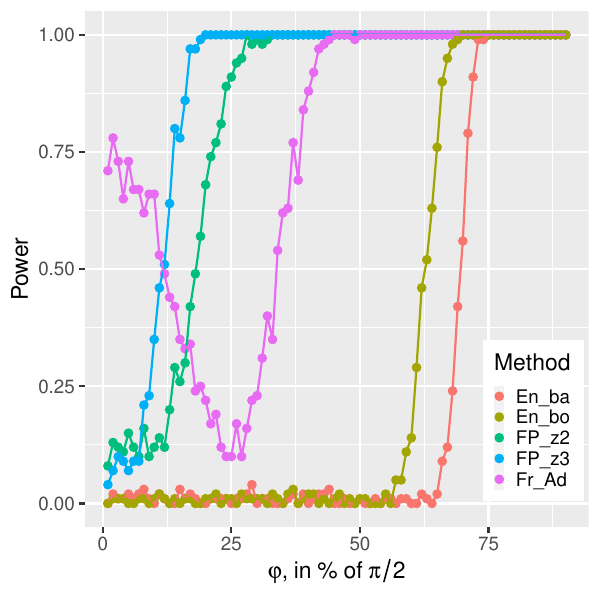}
			\label{t26}
		\end{subfigure}
		\begin{subfigure}{0.48\linewidth}
			\centering
			\caption{n=500,N=1000}
			\includegraphics[width=0.96\linewidth]{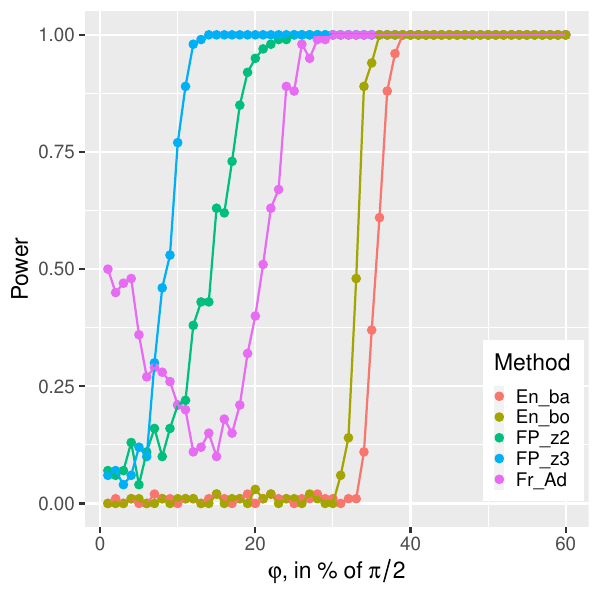}
			\label{t27}
		\end{subfigure}
		\centering
		\caption{Power comparison for Scenario II when $r_n=11$.}
		\label{tt13}
	\end{figure}

Figures \ref{tt12} and \ref{tt13} illustrate the power comparison in Scenario II when $r_n=7$ and $r_n=11$, respectively. Similar to Scenario I, our statistics maintain satisfactory performance, and both En$\underline{~}$Bo and En$\underline{~}$Ba show significant power loss especially when $\varphi$ is small, say less than $\pi/2\times0.4$. One may notice that Fr$\underline{~}$Ad demonstrates the highest power under an extremely weak alternative (e.g., $\varphi=\pi/2\times0.01$). However, we mention
	that this high power may not be trusted due to its empirical size being much larger than the nominal level $0.1$ as observed from Table \ref{table2}.

\begin{figure}[H]		
		\centering
		\begin{subfigure}{0.48\linewidth}
			\centering
			\caption{n=500,N=500}
			\includegraphics[width=0.96\linewidth]{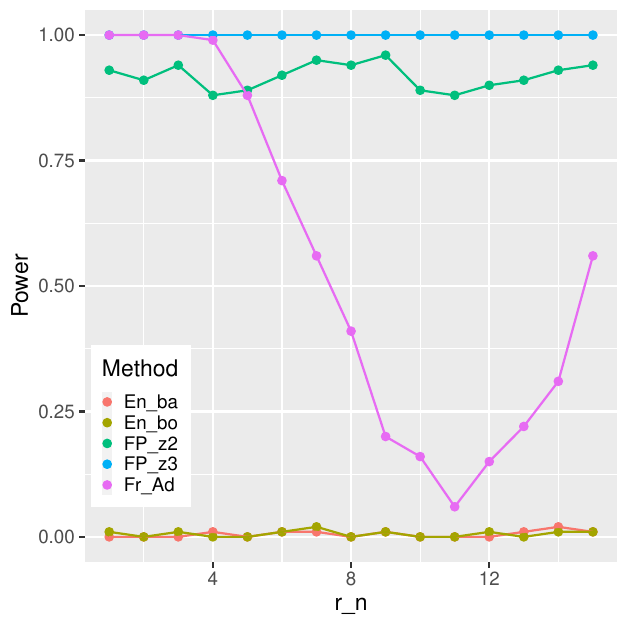}
			\label{t28}
		\end{subfigure}
		\begin{subfigure}{0.48\linewidth}
			\centering
			\caption{n=500,N=1000}
			\includegraphics[width=0.96\linewidth]{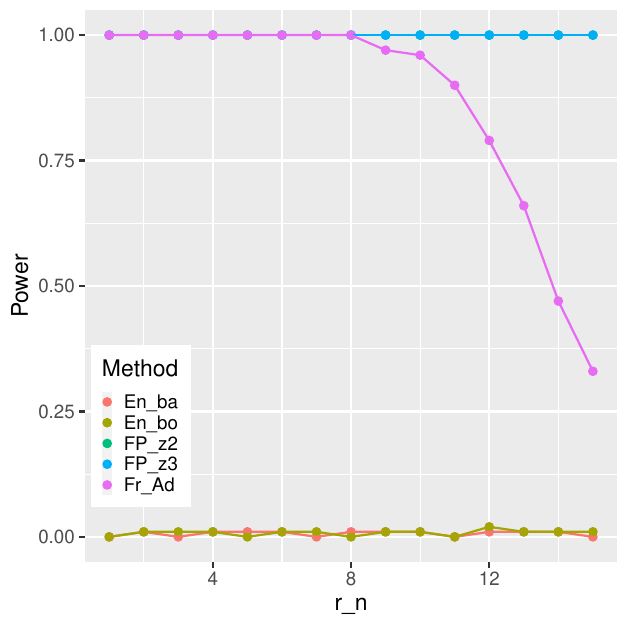}
			\label{t29}
		\end{subfigure}
		\centering
		\caption{Power comparison for Scenario III when the
			angle $\varphi=\pi/8$. The rank $r_n$ varies within $\{1,2,\cdots,15\}$.}
		\label{tt14}
	\end{figure}

\begin{figure}[H]		
		\centering
		\begin{subfigure}{0.48\linewidth}
			\centering
			\caption{n=500,N=500}
			\includegraphics[width=0.96\linewidth]{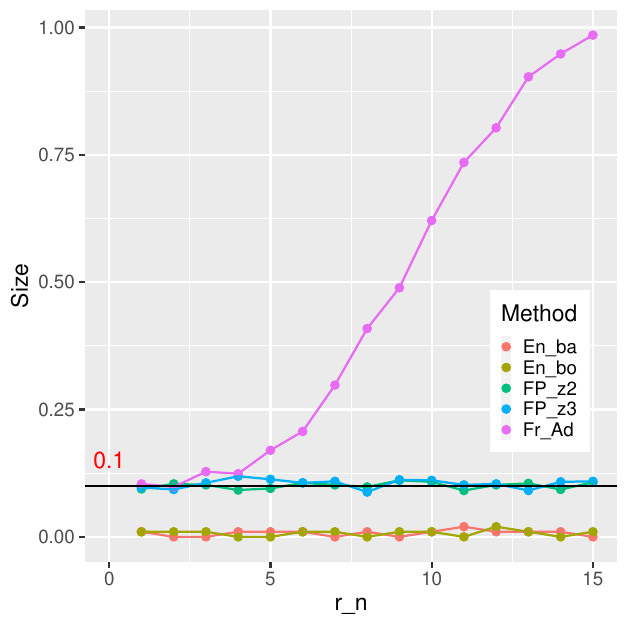}
			\label{t30}
		\end{subfigure}
		\begin{subfigure}{0.48\linewidth}
			\centering
			\caption{n=500,N=1000}
			\includegraphics[width=0.96\linewidth]{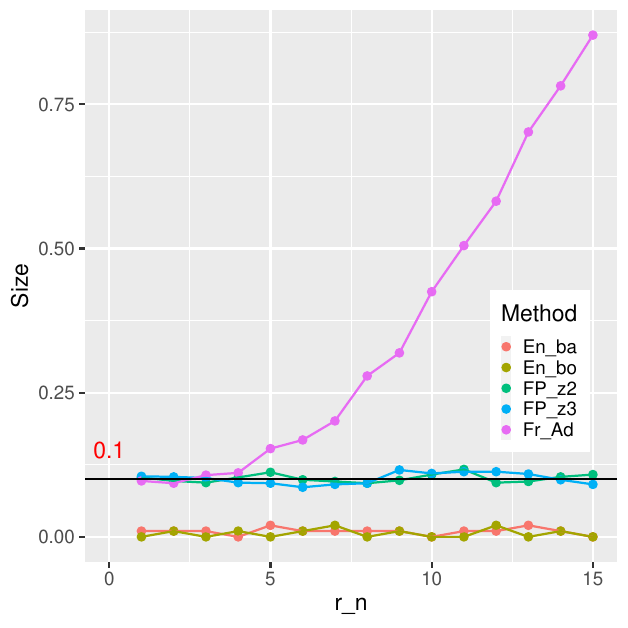}
			\label{t31}
		\end{subfigure}
		\centering
		\caption{Empirical sizes for Scenario III when the angle $\varphi=0$. The rank $r_n$ varies within $\{1,2,\cdots,15\}$. The black line is the nominal level $\alpha=0.1$.}
		\label{tt15}
	\end{figure}

Figure \ref{tt14} displays the the power performances of these methods when $\varphi=\pi/8$ in Scenario III. Our statistic FP$\underline{~}{z^3}$ demonstrates superior power performance, especially for large rank $r_n$. We observe that the power of Fr$\underline{~}$Ad exhibits excellent performance for small values of rank $r_n$, but experiences a significant decline as $r_n$ increases. The powers of both En$\underline{~}$Bo and En$\underline{~}$Ba are close to zero across all $r_n$.
The empirical sizes corresponding to Scenario III when $\varphi=0$ are depicted in Figure \ref{tt15}. It is evident that both our methods FP$\underline{~}{z^3}$ and FP$\underline{~}{z^2}$ consistently exhibit accurate distribution, with empirical sizes closely approximating 0.1.
Fr$\underline{~}$Ad experiences inflated sizes as $r_n$ increases, while the sizes of En$\underline{~}$Bo and En$\underline{~}$Ba remain close to zero.
	
	\begin{table}[htb]
          \setlength{\tabcolsep}{1.3mm}
			\begin{center}
			\caption{Empirical sizes at the nominal level $\alpha=0.1$, based on 1000 replications. The two values closest to 0.1 are highlighted in bold.}
			\label{table2}
			\renewcommand\arraystretch{1.7}
		\begin{tabular}{|c|c|c|c|c|c|c|c|c|c|c|}
				\hline & \multicolumn{5}{c|}{$N=500$} & \multicolumn{5}{c|}{$N=1000$} \\
				\hline Method & FP-$z^2$ & FP-$z^3$ & Fr-Ad & En-Bo & En-Ba & FP-$z^2$ & FP-$z^3$ & Fr-Ad & En-Bo & En-Ba\\
				\hline Scenario I: $\mathcal{N}(0,1)$ &\textbf{ 0.105} & \textbf{0.094} & 0.108 & 0.013 & 0.005 & \textbf{0.094} & \textbf{0.095} & 0.109 & 0.006 & 0.007 \\
				\hline Scenario I: $t(10)$ & 0.104 & \textbf{0.103} & \textbf{0.099} & 0.010 & 0.009 & \textbf{0.096} & \textbf{0.097} & 0.112 & 0.003 & 0.005 \\
				\hline Scenario II: $r_n=7$ & \textbf{0.095} & \textbf{0.104} & 0.271 & 0.009 & 0.008 & \textbf{0.096} &\textbf{0.102} & 0.216 & 0.010 & 0.004 \\
				\hline Scenario II: $r_n=11$ & \textbf{0.101} & \textbf{0.091} & 0.720 & 0.006 & 0.006 & \textbf{0.103} & \textbf{0.095} & 0.508 & 0.007 & 0.009 \\
				\hline
			\end{tabular}
		\end{center}
	\end{table}	

\subsection{An empirical examination on the universality of the functional projection approach}\label{addsimu}
In this section, we conduct an empirical examination to demonstrate the effectiveness of our functional projection approach (\ref{302}) when all $d_i$'s are smaller than $\sqrt{c_n}$, as mentioned in Remark \ref{rmk42}. As in Section \ref{sub5.1}, we assume that under the null hypothesis $\boldsymbol{H}_0$, the eigenvectors align with the axes of the coordinate system, i.e., $\bm{v}_i=\bm{e}_i$ for $i=1,\cdots,r_n$. The
hypothetical projection matrix and default covariance matrix are given in (\ref{HPM}) and (\ref{DCM}). For the alternative, we rotate the first $r_n$ eigenvectors by an angle $\varphi$. To be more specific, the covariance matrix under $\bm{H}_1$ can {be} explicitly written as
\begin{align*}
    \bm{\Sigma}_n=\bm{I}_n+\sum_{i=1}^{r_n}d_i\bm{v}_i^{\varphi}(\bm{v}_i^{\varphi})^{\top},
\end{align*}
where
\begin{align*}
    \bm{v}_i^{\varphi}=(\underbrace{0,\cdots,0}_{i-1},\underbrace{\cos{\varphi},0,\cdots,0}_{r_n},\sin{\varphi},0,\cdots,0)^{\top}.
\end{align*}
Consider the following scenario:

\textbf{$\bullet$ Scenario \uppercase\expandafter{\romannumeral4}.}
	Set $d_1=d_2=\cdots=d_{r_n}=0.5$.  $X_{ij}\sim \mathcal{N}(0,1)$ is taken into account. The dimension $n=500$ and the sample size $N=1000$. Obviously, all $d_i<\sqrt{c_n}$. Under the null hypothesis, the rank $r_n$ varies within $\{1,\cdots,15\}$ to check the distribution accuracy. Under the alternative hypothesis, we consider $r_n=5, 7, 9$ and vary the angle $\varphi$ within $\{1\%, 2\%, \cdots, 100\%\}\times \pi/2$ to capture  power performance trends.

\begin{figure}[H]		
		\centering
		\begin{subfigure}{0.48\linewidth}
			\centering
			\caption{Empirical power}
			\includegraphics[width=0.96\linewidth]{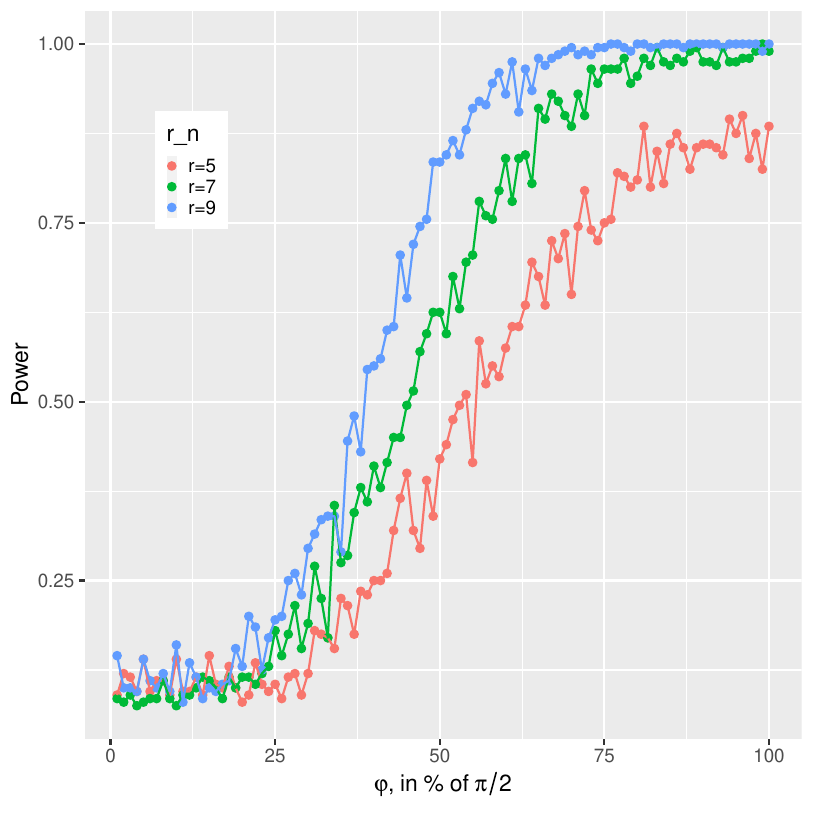}
			\label{t200}
		\end{subfigure}
		\begin{subfigure}{0.48\linewidth}
			\centering
			\caption{Empirical size}
			\includegraphics[width=0.96\linewidth]{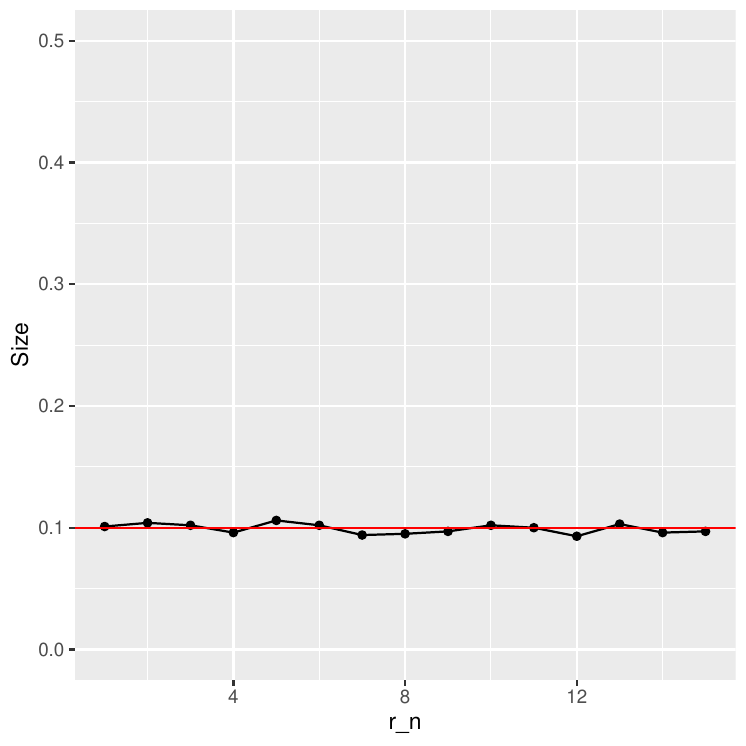}
			\label{t201}
		\end{subfigure}
		\centering
		\caption{Empirical powers (left panel) and sizes (right panel) for Scenario \uppercase\expandafter{\romannumeral4} when all $d_i<\sqrt{c_n}$. The red line in plot (b) is the nominal level $\alpha=0.1$.}
		\label{ttnew}
	\end{figure}

The choices for the remaining parameters are as follows: the nominal level $\alpha=0.1$, the threshold $\delta$ in (\ref{esti}) is $\delta=0.1$, and the function $f(z)=z^3$. Empirical powers are calculated from $200$ replications, while empirical sizes are recorded based on $1000$ replications. The results presented in Figure \ref{ttnew} demonstrate a strong alignment between our function projection approach and the theoretical normal distribution under the null hypothesis, with varying values of $r_n$. Furthermore, our method exhibits enhanced power performance as $r_n$ and $\varphi$ increase.

\begin{acks}[Acknowledgments]
	Yanlin Hu and Qing Yang are co-first authors. Xiao Han is the corresponding author. The authors would like to thank the anonymous referees, an Associate
	Editor and the Editor for their constructive comments that improved the
	quality of this paper. 
\end{acks}

\begin{funding}
	This work was supported by National Natural Science Foundation of China (Grant No. 12571297), National Natural Science Foundation of China (Grant No.12371278), National Key R\&D Program of China-2022YFA1008000 and the Talents Introduction Program of the Chinese Academy of Sciences (Category B).
\end{funding}

\begin{supplement}
	\stitle{Supplementary material for  ``Generalized Linear Spectral Statistics of High-dimensional Sample Covariance Matrices and Its Applications''}
	\sdescription{The supplementary material contains additional results on simulation results and all the technical proofs.}
\end{supplement}

\bibliographystyle{imsart-number}
	\bibliography{references}	

\end{document}